\documentclass[10pt,a4paper]{article}
\usepackage[utf8]{inputenc}
\usepackage[english]{babel}
\usepackage{amsmath}
\usepackage{amsfonts}
\usepackage{amssymb}
\usepackage[english]{babel}
\usepackage{color}
\usepackage{array}
\usepackage{enumerate}
\usepackage{enumitem}
\usepackage{refcount}
\usepackage{algorithm}
\usepackage{algorithmic}
\usepackage{boxedminipage}
\usepackage{float}
\usepackage{graphicx}
\usepackage{multirow}
\usepackage{tikz}
\usepackage{amssymb}
\usepackage{caption}

\newtheorem{thm}{Theorem}[section]

\newtheorem{lem}[thm]{Lemma}

\newtheorem{defn}[thm]{Definition}

\newtheorem{exa}[thm]{Example}

\newtheorem{remark}[thm]{Remark}

\newtheorem{NP cases}[thm]{Possible cases of Newton Polygons for a triple point with a triple tangent}

\newtheorem{stp}[thm]{Stop criterion}

\def\proof{\b {\bf Proof} }

\newcommand{\bC}{\mathcal C} %bC:=base C

\newcommand{\bR}{\mathcal R} % riferimenti affini

\newcommand{\Mod}[1]{\ (\mathrm{mod}\ #1)}

\newcommand{\prfend}{\hbox to7pt{\hfil}
\par\vskip-\baselineskip\hbox to\hsize
{\hfil\vbox {\hrule width6pt height6pt}}\vskip\baselineskip}

\def\a{\bigskip \par \noindent}
\def\b{\par \noindent}

\def\dd{\medskip \par \noindent}
\long\def\eatit#1{}

\def\Q{\Bbb Q}
\def\Z{\Bbb Z}
\def\N{\Bbb N}
\def\C{\Bbb C}
\def\A{\Bbb A}

\def\g{\gamma}
\def\G{\Gamma}
\font\tengothic=eufm10
\font\sevengothic=eufm7
\newfam\gothicfam
       \textfont\gothicfam=\tengothic
       \scriptfont\gothicfam=\sevengothic

    \font\tenmsb=msbm10              \font\sevenmsb=msbm7
\newfam\msbfam
       \textfont\msbfam=\tenmsb
       \scriptfont\msbfam=\sevenmsb
\def\Bbb#1{{\fam\msbfam #1}}

\usepackage{xypic}

\begin{document}
	\title{ The Newton-Puiseux Algorithm and triple points for plane curves}
\author{Stefano Canino, Alessandro Gimigliano,
	Monica Id\`a\\
}

\maketitle
\begin{abstract}
\noindent The paper is an introduction to the use of the classical Newton-Puiseux procedure, oriented to an algorithmic description of it. This procedure enables to get polynomial approximations for parameterizations of branches of an algebraic plane curve at a singular point. We look for an approach that can be easily grasped and almost self contained.  We illustrate the use of the algorithm, first in a completely worked out example of a curve with a point of multiplicity 6, and secondly in the study of triple points on reduced plane curves. \end{abstract}



%%%%%%%%%%%%%%%%%%%% HO fatto eatit{} sui tripli 
\section{Introduction}
%The aim of this paper is to study the structure of triple points for an algebraic plane curve: for such singularities one does not have standard analytic models for any possible case, as it happens for double points. 
\dd The study of singular points of algebraic plane curves is a very classical subject for Algebraic Geometry, and a renewed interest for this subject comes by the connections with applications (e.g. in CAD), in particular for rational curves, which (since they can be parameterized) are particularly suited for plotting a design. So to study ways to describe singular points of plane rational curves, and to find algorithms that can do that, have been quite a subject of research, and many results have been obtained, by attacking the problem from different points of view (e.g. see \cite{CWL}, \cite{CKPU} and \cite{BGI}).   
\dd When the curve is not rational, one does not have the advantage of a parameterization, nevertheless there are local parameterizations of the branches of the curve at a singular point, not in terms of polynomials, but by using power series. Namely, if the curve is defined by a polynomial $f(x,y) \in {\mathbb C}[x,y]$, with $f(0,0)=0$, one can find parameterizations  $T \rightarrow (T^m, p(T))$, such that $f(T^m,p(T))=0$, $p(0)=0$ and $p(T)\in {\mathbb C}\{T\}$, the ring of convergent power series in $T$. Every such parameterization describes a branch of $f(x,y)$ at $(0,0)$. 
\dd The way to get these parameterizations is via Puiseux series; these are quite a classical subject (the idea goes back to Newton, then rediscovered by Puiseux) and sometimes is taught even in undergraduate classes, but the choice to make it part of a regular course of study is not so common and is quite scattered, hence their knowledge is not so widespread among graduate students and young researchers (in Algebraic Geometry or connected studies).

Even the definition of what a Puiseux series is is not standard: we will consider them as power series in which also positive rational exponents appears and we denote with $\mathbb{C}\{\{x\}\}$ the ring of such series (that is our choice: even for this ring there is not a standard symbol). In some text the word {\it Puiseux series} is used for series where also a finite number of negative rational exponents is allowed and in this case Puiseux series would form a {\it field} (the field of quotients of our ring $\mathbb{C}\{\{x\}\}$) and the main Theorem in this theory (also called Newton-Puiseux Theorem) can be stated as {\it The field of Puiseux series is algebraically closed} (e.g. see \cite{W}). 
\dd Detailed treatises of the Theory (apart from the early reference of Newton and Puiseux themselves, see \cite{Ne} and \cite{P}) can be found for example in \cite{Ar} chs. 6,7, \cite{Fi}, chs. 7,8 and appx. 4, \cite{GLS1}, ch. 1 sec. 3, \cite{W} ch. IV., and it is syntethized in the {\it Newton-Puiseux Algorithm} (e.g. see either \cite{W} for a more theoretical approach or \cite{Du} , \cite{PR} and references therein for efficient computational versions).
\dd The aim of this paper is to give a handy version of the Newton-Puiseux procedure which avoids involving negative exponents and it is written as an ``as concise as possible"  Algorithm, together with a self-contained justification of why it actually works (although we use \cite{W} as main reference). This will be done in section 3, after a section of preliminaries (section 2). The Algorithm gives polynomial approximations for the branches of a reduced algebraic curve at a singular point.  We also give an algorithm to be used in case one knows the irreducible decomposition of the curve, and that can be used also when $f$  is not reduced. 
\dd We illustrate as the given Algorithm works, first with an example of a curve with several branches  at a point of multiplicity 6 (at the end of section 3) and then, in section 4,  by  giving the description of all possible structures of a triple point  with a triple tangent (see Theorem 4.3). When $C$ has only a 3-branch (the most interesting case) we see that its parameterization can be characterized by its type $s$ ($s\geq 4$ and $s\not\equiv 0$ (mod 3): see Definition \ref{type}), where the branch is analytically equivalent to $(T^3,a_sT^s+ \dots)$. 
\section{Preliminaries}
 In the following $\bC:f=0$ is a reduced plane algebraic curve in $\A^2(\C)$ and $O:=(0,0)$ is a point of $\bC$. Given a formal power series $s(T)\in\C [[T]]$, $o(s)$ denotes its order.
\begin{defn}\rm
A {\em formal parameterization of $\bC$ at $O$} is a pair $(p(T),q(T))$ with $p(T),q(T)\in\C [[T]] $ such that $f(p(T),q(T))=0,p(0)=0$ and $q(0)=0$.
\end{defn}
An equivalence relation can be defined on the set of formal parameterizations saying that two formal parameterizations $(p_1,q_1)$ and $(p_2,q_2)$ are  \textit{equivalent} if there exists an invertible series $s\in\C [[T]] $ such that $(p_1,q_1)=(p_2\circ s,q_2\circ s)$. Moreover, a formal parameterization is called \textit{reduced} if it is not equivalent to a formal parameterization of the form $(p(T^m),q(T^m)),m>1$.
\begin{defn}\rm
A {\em branch of $\bC$ at $O$} is a class of equivalent formal reduced parameterizations of $\bC$ at $O$.
\end{defn}
\begin{remark}\label{equiv}\rm
It can be shown (see for example \cite{W} ch. IV, Theorem 2.2) that in a suitable coordinate system any given parameterization is equivalent to one of the type
$$(T^n,p(T))$$
where $p(T)\in\C [[T]]$. Moreover, two such parameterizations 
$$(T^n,p(T)), \quad (T^m,q(T))$$
 are equivalent if and only if $n=m$ and there exists an $n$-th root of the unity $\omega$ such that $p(\omega T)=q(T)$. 
\end{remark}
\begin{remark}\rm
Sometimes, what we defined as a branch of $\bC$ at $O$ is called a {\it place} of $\bC$ at $O$ and can be viewed as the algebraic counterpart of the more usual analytical definition of branch. However, the two definitions are equivalent; more details can be found in \cite{W} pp. 96,97.
\end{remark}
\begin{defn}\label{tangent}\rm
Let $R=(p(T),q(T))$ be a branch of $\bC$ at $O$. We define the {\em multiplicity of $R$ at $O$} to be the number $m=m_O(R):=\min\left\lbrace o(p),o(q)\right\rbrace $ and we say that $R$ is a $m$-branch of $\bC$. Moreover, if $c$ and $d$ are the coefficients of $T^{m}$ respectively in $p$ and $q$, the line passing through $O$ with direction $(c,d)$ is called the tangent line of $R$.
\end{defn}
\begin{remark}\label{mult}\rm
It can be checked that the equation of the tangent cone to $\bC$ at $O$ is the product of the tangent lines (possibly repeated) to the different branches of $\bC$ at $O$, each one raised to the multiplicity of the corresponding branch. It follows that:
\begin{itemize}
	\item $\bC$ has in each of its points a finite number of branches.
	\item  If $t$ is a tangent of $\bC$ at $O$ with multiplicity $m_t$ and $R_1,\dots,R_s$ are the branches of $\bC$ at $O$ whose tangent is $t$ then $\sum_{i=1}^{s}m_O(R_i)=m_t$
	\item If $R_1,\dots,R_n$ are the branches of $\bC$ at $O$ then $\sum_{i=1}^{n}m_O(R_i)=m_O(\bC)$.
\end{itemize}
\end{remark}
In order to find explictly the branches of a curve, we will use the Newton-Puiseux algorithm, which uses Puiseux series and the Newton polygon. The former are nothing but power series in which positive rational exponents are admitted. More formally we have
\begin{defn}\rm
A \textit {Puiseux series} is an element $p \in \C\left\lbrace \left\lbrace x\right\rbrace \right\rbrace: =\bigcup_{r=1}^{\infty}\C[[ x^{\frac{1}{r}}]]  $. The minimum $r$ such that $p\in\C[[ x^{\frac{1}{r}}]] $ is called {\it the lowest common denominator of }$p\in\C\left\lbrace \left\lbrace x\right\rbrace \right\rbrace $ (it is also called  order of polydromy, e.g. see \cite{Ma}). The {\it  order of $p$} is defined analogously to the formal series order.
\end{defn}
Notice that a finite sum of Puiseux series is a Puiseux series, since if $s_1\in \C [[x^{1\over {r_1}}]], \dots, s_n\in \C [[x^{1\over {r_n}}]]$, then $s_1+\dots + s_n \in  \C [[x^{1\over {r_1}\dots r_n}]]$.
\begin{defn}\rm
A \textit {Puiseux $y$-polynomial} is a polynomial  $f\in\C\{\{x\}\} [y]$ of the form
$$f(x,y)=\sum_{i=0}^{n}p_i(x)y^i$$
where each $p_i(x)$ has a finite number of terms,  i.e. $p_1\in \C [x^{1\over {r_1}}], \dots, p_n\in \C [x^{1\over {r_n}}]$.
Setting $r=r_1\dots r_n$ we have that $p_1,\dots, p_n \in  \C [x^{1\over {r}}]$, hence $f$ is polynomial in $\C[x^{1\over {r}},y]$. 
\end{defn}
A root of $f$ is a series $p(x^{1\over q})$ such that $f(x, p(x^{1\over q}))=0$.
\begin{defn}\rm
A Puiseux $y$-polynomial $f\in\C[x^{1\over r}, y]$ is called \textit{Newton convenient} if $x^{1\over r}$ and $y$ do not divide $f$ (i.e. there are $a,b\in\N$ such that the monomials $x^{a\over r}$ and $y^b$ appear in $f$ with non zero coefficients), and $f(O)=0$. 
\end{defn}
The Newton polygon can actually be defined for polynomials in any number of variables and, in general, is called \textit{Newton diagram} (see \cite{GLS1} Definition 2.14 pp. 121-122). We use only the \linebreak 2-dimensional Newton diagram of Newton convenient Puiseux $y$-polynomials and we refer to it as \textit{Newton polygon} (see also \cite{W}).
\begin{defn}\label{NPoly} \rm
	Let $\displaystyle f=\sum_{i+j=1}^{d}a_{ij}(x^{1 \over r})^iy^j\in\C[x^{1 \over r},y]$ be a Newton convenient Puiseux $y$-polynomial and $\text{supp}(f):=\{({i\over r},j)\in\Q^2,a_{ij}\neq0 \}$ the support of $f$. We construct the {\em Newton Polygon of $f$}, denoted by $\Gamma(f)$, as follows:
	\begin{itemize}
		\item Let $P_1$ be the leftmost point of $\text{supp}(f)$  (notice that $P_1$ is on the $y$-axis; if there are more points with abscissa 0 we  choose the lowest one). From $P_1$ we rotate a vertical downward ray counterclockwise and stop rotating it when it meets the first point of $\text{supp}(f)$. We call $P_2$ the rightmost point of $\text{supp}(f)$ met by this ray. 		
		\item From $P_2$ we rotate a vertical downward ray counterclockwise and we stop rotating it when it meets the first point of $\text{supp}(f)$. We call $P_3$ the rightmost point of $\text{supp}(f)$ met by this ray.
		\item We repeat the procedure until we reach a point  $P_{k+1}$ on the $x$-axis.
			\end{itemize}
 For all  $P_n,P_{n+1}$, we denote by $\textbf{a}_n$ the segment $\overline{P_nP_{n+1}}$, we set $\Gamma(f)=\{\textbf{a}_1\dots\textbf{a}_k\}$ and we say that $\textbf{a}_1,\dots,\textbf{a}_{k}$ are the edges of $\Gamma(f)$. The difference between the ordinates of $P_n$ and $P_{n+1}$ is called {\em height of} $\textbf{a}_n$. For each $\textbf{a}\in\Gamma(f)$ we define the \textbf{a}-truncation of $f$ as 
		$$ f^{\textbf{a}}=\sum_{({i\over r},j)\in\textbf{a}\cap\Q^2}a_{ij}(x^{1\over r})^iy^j.$$
	\end{defn}
 \begin{remark}\label{novh} \rm In the notations of the previous definition, consider the unbounded subset $N$ of the first quadrant $x\geq 0, y \geq 0$ whose border is given by the Newton Polygon union with the ray of the $y$-axis of vertex $P_1$ and the ray of the $x$-axis of vertex $P_{k+1}$. By construction $N$ is convex; since the Newton Polygon starts at the point $P_1$ of the $y$-axis and ends at the point $P_{k+1}$ of the $x$-axis, it's clear that there are no vertical neither horizontal edges in the Newton Polygon, and that all the edges have negative slope.
\end{remark} 
The next well known theorem will be very useful in the following (its proof can be found in \cite{Ar}, Theorem 6.4).
\begin{thm}[Implicit Function Theorem for Puiseux series]\label{ift} Let $$f=\sum_{i=0}^dp_i(x)y^i\in\C\left\lbrace \left\lbrace x\right\rbrace \right\rbrace [y].$$
If $p_0(0)=0$ and $p_1(0)\neq0$ then $\exists!\;p\in\C\left\lbrace \left\lbrace x\right\rbrace \right\rbrace$  such that $p(0)=0$ and $f(x,p(x))=0$. Moreover, if $p_i\in\C [[x^{\frac{1}{r}}]]$ for $i=1,\dots, d$, then $p\in\C [[ x^{\frac{1}{r}}]] $. 
\end{thm}
\section{The Newton-Puiseux algorithm}
Now we want to give the algorithm (Newton-Puiseux Algorithm) which is the aim of this section: it will determine, given a reduced polynomial $f=\sum_{i+j=1}^da_{ij}x^iy^j$, all the branches of the algebraic plane curve $\bC: f=0$ at $O$. With this in mind,  we give some preliminaries needed for the algorithm.
 \subsection{The $\ast $-procedure and the graph $\textbf{G}_f$}\label{NewP}
In the following  $h$ denotes a Puiseux $y$-polynomial with $h(O)=0$, $h\in\C[x^{1\over r}, y]$, and such that $x^{1\over r} \nmid h$  (recall that $r\in \N, r>0$). Hence we consider powers of $x$ with rational exponents, while the powers of $y$ always have an integer exponent.
\dd 
At this moment we just describe the procedure and its successive applications, later on we shall give the needed justifications. In the procedure we distinguish two cases, according to whether  $y\mid h$ or not.
	
\a { $\underline{\bf {Case}\quad y\nmid h:}$  By assumption $x^{1\over r} \nmid h$, $h(O)=0$ and $y\nmid h$, hence $h$ is Newton convenient.
	We consider the Newton Polygon $\Gamma(h)=\left\lbrace \textbf{b}_1,\dots, \textbf{b}_k\right\rbrace $ of $h$. We set
	$$h^{(i)}:=\frac{h^{\textbf{b}_i}}{x^{u_i}y^{v_i}} \qquad \forall\; i=1,\dots, k$$
	where $x^{u_i}$ and $y^{v_i}$ are respectively the highest power of $x$ and $y$ such that $h^{(i)}\in\C[x^{1\over r}, y]$. Moreover, we have $h^{(i)}(O)=0$ {(this will be justified). 
		\dd We find all the $s_i$ roots $q_{i1},\dots q_{is_i}$ of $h^{(i)}$, considered as a polynomial in $y$;  if $-\frac{1}{r_i}$ is the slope of ${\textbf{b}_i}$, then $h^{(i)}$ has all its roots in $\C[x^{ r_i}]$ and each $q_{ij}$ is of the form $q_{ij}(x)=c_{ij}x^{r_i}$ with $c_{ij}\in \C$ (this will be justified). 
	\b We finally set
	$$h^{(ij)}:=\frac{h(x,x^{r_i}(c_{ij}+z))}{x^{m_{i}}}\qquad \forall \; i=1,\dots, k,\quad \forall\;  j=1,\dots,s_i$$
	where $m_{i}$ is the highest power of $x$ such that $h^{(ij)}\in\C\left\lbrace \left\lbrace x\right\rbrace \right\rbrace [z]$.
	\dd   By construction, $h^{(ij)}$ is a Puiseux $z$-polynomial, and if $h^{(ij)}$ is in, say, $\C[x^{1\over b}, y]$, then $x^{1\over b} \nmid h$. Moreover, we have $h^{(ij)}(O)=0$ (this will be justified). So the result of the procedure is:
	$$s_i= s_i(h),\quad q_{ij}= q_{ij}(h), \quad c_{ij}= c_{ij}(h), \quad r_i= r_i(h), \quad h^{(ij)}, \quad i=1,\dots, k,\quad  j=1,\dots,s_i.$$
	
\a { $\underline {\bf{Case}\quad y\mid h}$: Let $y^e$ be the highest power of $y$ such that $\hat{h}:=\frac{h(x,y)}{y^e}\in\C\left\lbrace \left\lbrace x\right\rbrace \right\rbrace [y]$. 
 We have  two subcases: $\hat{h}(O)\neq 0$ and $\hat{h}(O)= 0$.
\begin{itemize} 
\item If $\hat{h}(O)\neq 0$ we ignore $\hat h$ and we take into account uniquely the contribution given by the root $y=0$ of $h$, in the following way: we define a virtual Newton Polygon $\Gamma(h)$ with just one edge, so we set, with an abuse of notation:  
	 $$ \Gamma(h):= \{  \textbf{b}_1 \}, \quad s_1:=1, \quad q_{11}:=0, \quad c_{11}=0, \quad r_{1}:=0, \quad  h^{(11)}:=0.$$
\item If $\hat{h}(O)=0$ we take into account the contribution given by $\hat{h}$ as well as the one given by the root $y=0$, in the following way: 
	$\hat h$ being Newton convenient, we can consider the Newton Polygon $\Gamma(\hat h)=\left\lbrace \textbf{b}_1\dots\textbf{b}_k\right\rbrace $ of $\hat h$
	 and apply the procedure defined in Case $y\nmid h$ to the Puiseux $y$-polynomial $\hat h$; we rename the results of the procedure as follows:  
	$$s_i(h):= s_i(\hat h),\; q_{ij}(h)= q_{ij}(\hat h), \;  c_{ij}(h):= c_{ij}(\hat h), \;  r_i(h):= r_i(\hat h), \;  h^{(ij)}:=\hat h^{(ij)}, \;  i=1,\dots, k,\;   j=1,\dots,s_i.$$
	Moreover we add a further virtual edge $\textbf{b}_{k+1}$ and we set, with an abuse of notation: 
	$$\Gamma(h):= \{  \textbf{b}_1, \dots,  \textbf{b}_k,\textbf{b}_{k+1}\}, \quad s_{k+1}:=1, \quad q_{(k+1)1}:=0, \quad c_{(k+1)1}=0, \quad r_{k+1}:=0, \quad h^{((k+1)1)}:=0.$$
	\end{itemize}
	
	\medskip
 Summarizing, the input of the $\ast $-procedure is a Puiseux $y$-polynomial $h\in\C[x^{1\over r}, y]$ such that $x^{1\over r} \nmid h$ and $h(O)=0$, and the output is:
$$s_i= s_i(h)\in \N,\quad r_i=r_i(h)\in \Q_+,\quad c_{ij}=c_{ij}(h)\in \C$$ and some Puiseux $z$-polynomials $h^{(ij)}\in\C\left\lbrace \left\lbrace x\right\rbrace \right\rbrace [z]$, for any choice of the edge  $\textbf{b}_i$ of $\Gamma(h)$ and for any $q_{ij}$ relative to that edge. 

Moreover, each polynomial $h^{(ij)}$ in the output either satisfies the necessary conditions to apply again the \linebreak $\ast $-procedure, or $h^{(ij)}=0$.
\dd
Now assume a polynomial $f\in \C[x,y]$ is given, with $f(O)=0$ and $x \nmid f$. We apply the $\ast$-procedure successively to all the polynomials $f,f^{(ij)},f^{(ij)\,(lt)},\dots$ getting a graph $\textbf{G}_f$ of the form

\bigskip
	\begin{figure}[H]
		\centering
		\includegraphics[scale=0.25]{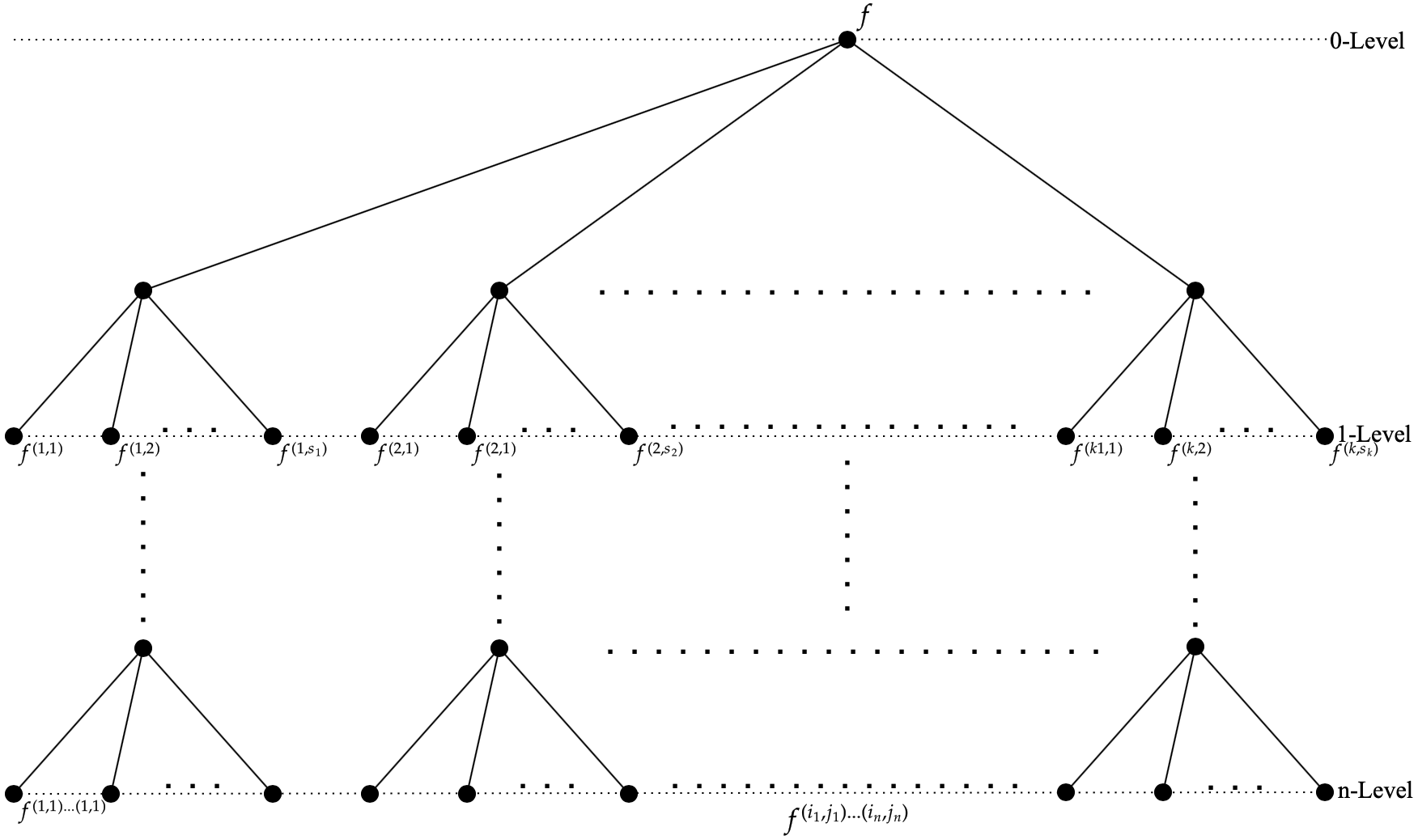}
		\caption {The graph $\textbf{G}_f$}
	\end{figure}
\noindent The 0-level consists only of the vertex $f$, the 1-level consists of the vertices $f^{(ij)}$ and so on, so that the $n$-level is composed by the vertices marked with $2n$ indices.
\dd Since, in general, Puiseux series have infinite terms, the $\ast$-procedure, in general, does not end without a stop criterion. 
\begin{stp}\label{stop2}	In order to have that $\textbf{G}_f$ is a finite graph, we stop the procedure every time that we meet a polynomial ${\overline{f}}$ such that one of the following conditions occurs:
\begin{enumerate}[label=(\roman*)]
	\item ${\overline{f}}=0$; 
	\item ${\overline{f}}$ satysfies the assumptions of Theorem \ref{ift}, i.e. the monomial $z$ appears in ${\overline{f}}$ with a non-zero coefficient (graphically the Newton Polygon $\Gamma({\overline{f}})$ has an unique edge of height one).
\end{enumerate}
\end{stp}
 \subsection{The paths on $\textbf{G}_f$ and the branches of  the curve}\label{paths} 	
We consider now each possible descending path $\g$ in $\textbf{G}_f$ starting from $f_0:=f$ and we stop when we reach an $f^{(i_1\,j_1)\dots(i_n\,j_n)}$ which satisfies $(i)$ or $(ii)$ of \ref{stop2} above. Notice that we get in general paths of different lengths.

\medskip
 For each path $\g$, we call $n$-th step the passage from the $n$-level to the $(n+1)$-level of $\g$; graphically it is given by 3 vertices and 2 edges: the highest vertex is a polynomial $f^{(i_1j_1)\dots (i_nj_n)}$ which from now on is denoted by $f_n(\g)$ and, since the $\ast$-procedure introduces a new variable at each step, we consider $f_n(\g)\in\C\left\lbrace \left\lbrace x\right\rbrace \right\rbrace [y_n]$. The middle vertex is associated to an edge  of $\Gamma(f_n)$ which is now called $\textbf{a}_n(\gamma)$ and the lowest one is now called $f_{n+1}(\g)$. When no confusion is possible we will just write $f_n,\textbf{a}_n,f_{n+1}$.
\dd Once a path is fixed, we use the following notations for the $\ast$-procedure at the $n$-th step, assuming that $f_n\in \C[x^{{1\over r}},y_n]$, where $r$ depends on $\g$ and $n$ and that $\textbf{a}_n$ is not a ``virtual edge" which sprouts from a zero root:
	\begin{itemize}
		\item $x^{u_n}$ and $y_n^{v_n}$ are respectively the highest power of $x$ and $y_n$  such that: 
		$$g_n=\frac{f_n^{\textbf{a}_n}(x,y_n)}{x^{u_n}y_n^{v_n}}\in\C\{\{x\}\}[y_n].$$
		\item $q_n=c_nx^{r_n}$ is the chosen root of $g_n$.
		\item $m_n$ is the highest power of $x$ such that 
		$$f_{n+1}=\frac{f_n(x,x^{r_n}(c_n+y_{n+1}))}{x^{m_n}}\in\C\{\{x\}\}[y_{n+1}].$$
	\end{itemize}
If, instead, $\textbf{a}_n$ is a ``virtual edge", we just have $c_n=0, \; r_n=0, \;  f_{n+1}=0$ (so the path $\g$ stops here by the stop criterion).
\a {\bf Each path $\g$ gives the approximation of a branch:}
	
	\dd Let us consider the algebraic curve $\bC: f(x,y)=0$ through the origin. In the notations above, the path $\g$ gives the approximation of a branch of $\bC$ at $O$  in the following way:
	\dd if $\overline{n}$ is the lenght of the path $\g$, we set 
	$$p_{\g}(x):=\sum_{n=0}^{\overline{n}}c_nx^{r_0+\dots+r_n}.$$
	If no confusion arises we just write $p(x)$ for $p_{\g}(x)$.  
	\dd Let $r$ be the common lowest denominator of $p(x)$; Theorem \ref{ift} guarantees that $p(x)$ is the truncation of a Puiseux series $\tilde{p}(x)$, with common lowest denominator $r$. Moreover, the way the algorithm is constructed guarantees that $f(x,\tilde p(x))=0$ and $\tilde p(0)=0$, so that  $(T^r,p(T^r))$ is the approximation of a branch $(T^r,\tilde{p}(T^r))$ of $\bC$ at $O$. For more details about this see \cite{W} pp. 98-99.
	 \dd The above procedure is summarized in the forthcoming Algorithm 1. Notice that in the Algorithm we start with a reduced curve $\bC: f(x,y)=0$, i.e. the polynomial $f=f_0$ is supposed to be without multiple factors. This, by Lemma \ref{lemma4}, implies that the case relative to a ``virtual edge" (which would occur when $f_{\overline{n}-1}$ has a zero root, forcing the path $\g$ to stop) can never occur, unless $\bC$ has a branch coming from a component that is a rational curve, i.e. there is an algebraic branch at $O$.	In other words, if $(i)$ of the Stop criterion \ref{stop2} occurs, it means that $y_{\overline{n}-1}\mid f_{\overline{n}-1}$, i.e. $0\in\C\left\lbrace \left\lbrace x\right\rbrace \right\rbrace $ is a root of $f_{\overline{n}-1}$: hence, we can stop following the path in question because the corresponding Puiseux series has a finite number of terms and we have found all of them.
	\dd Of course the roots of $f$ have, in general, infinite terms. Every time we apply the  $\ast$-procedure, we get a new polynomial which has a Newton Polygon whose height is $\leq$ the previous one by Lemma \ref{lemma2}. The criterion actually yields to a stop on any descending path in  $\textbf{G}_f$ because $f$ is reduced and this implies that each path on the graph ends with a step of height 1 (see also  \cite{W} p.105) hence we end in case $ii)$ of the step, unless we ended in case $i)$ before.
    \dd If, once arrived at the $\bar n^{th}$ step, we are in case $(ii)$ of the Stop criterion \ref{stop2}, i.e.  $f_{\overline{n}}(x,y_{\overline{n}})$ satisfies the hypothesis of Theorem \ref{ift}, we have that for any $n\geq\overline{n}$ there is a unique choice for $c_n$ and $r_n$. This means that the root $p(x)$ does not split anymore, i.e. if we were to go on along the path from now on there would be a unique possible choice at each step. Hence,  we have found enough terms of $p_\g(x)$ to distinguish it from the other branches, therefore we can stop the algorithm. Clearly, if one wants to calculate more terms of $ p(x)$ it is enough to continue the procedure.
 \dd In conclusion, the Stop criterion \ref{stop2} is a good one.
\a {\bf Justification of the procedure}

\a There are several results that are needed in order to show that we can actually perform the $\ast $-procedure described in Subsection \ref{NewP}, namely:
\begin{enumerate}[leftmargin=*,label=(\arabic*)]
	\item In the  $\ast $-procedure, we need $h^{(i)}(O)=0$, and the roots $q_{i1},\dots q_{is_i}$ of $h^{(i)}$ need to be in $\C[x^{1\over r}]$, with  $q_{ij}$ of the form $q_{ij}(x)=c_{ij}x^{r_i}$, $c_{ij}\in \C$;
	\item In order to keep applying the  $\ast $-procedure, we must show that the Puiseux $z$-polynomial $h^{(ij)}$ satisfies $h^{(ij)}(O)=0$. 
\end{enumerate}
Rephrasing $(1)$ and $(2)$ in the language used above in this Subsection \ref{paths}, we must show that:
\begin{enumerate}[leftmargin=*,label=(\arabic*)]
	\item $g_n(O)=0$, and the roots $q_{i1},\dots q_{is_i}$ of $g_n$ need to be in $\C[x^{1\over r}]$, with  $q_{ij}$ of the form $q_{ij}(x)=c_{ij}x^{r_i}$, $c_{ij}\in \C$.  
	\item $f_{n+1}(O)=0$.  
\end{enumerate}
Finally, notice that we need to distinguish when the parameterizations we get in the algorithm are equivalent, in which case they give the same branch. Lemma \ref{equivalent} gives a criterion to do that.	The fact that (1) and (2) above hold will be proved by Lemma \ref{lemma2} and Lemma \ref{lemma1} below, respectively.
 
\bigskip
 \subsection{The Algorithm}\label{just}
\a {\bf The Newton-Puiseux Algorithm}
\dd Now we give the Newton-Puiseux Algorithm based on the procedures described above. Notice that we suppose  $x\nmid f$, since we used this hypothesis in the $*$-procedure; of course if one wants to study a curve that has the $y$-axis as a component, it is sufficient to remove the factor $x$ from $f$ and then to add the vertical branch $\{x=0\}$ (parameterized by $(0,T)$) to the branches at $O$ obtained by the Algorithm.

\medskip
\noindent We give Algorithm 1 for reduced curves; notice that when $f$ is not irreducible, it can be much easier, computationally, to work separately on its irreducible components (if they are known: to find them is not computationally easy, e.g. see \cite{PW}, where an algorithm is given) and then consider the union of the branches of any such component.  We summarize this (quite obvious) procedure in Algorithm 2, where we also drop the ``$f$ reduced" hypothesis, by adding that any branch must be counted as ``$n$-ple" if obtained from a factor $f_i^n$.  Again we suppose that $x=0$ is not a component of $\bC$; obviously if  there is a factor of $f$ which is $x^n$, we just remove it before applying the Algorithm and then we add the  branch $\{x=0\}$ counted $n$ times to the output.

\begin{algorithm}[H]\caption{Study of the branches at $(0,0)$ of a reduced plane curve $\bC $ of equation $f=\sum_{i+j=1}^da_{ij}x^iy^j$, $a_{ij}\in \C$, such that $x\nmid f$}\label{Algorithm1}
\a  \textbf{Input}: $\bC $: $f=\sum_{i+j=1}^da_{ij}x^iy^j$, $O=(0,0)$.\\

 \textbf{Output}: Integers $r_1,...,r_s$ and $p_{1} \in \C [x^{1\over r_{1}}],...,p_s\in \C [x^{1\over r_{s}}]$  such that $\bC $ possesses $s$ branches at $O$, and each branch is approximated by $(T^{r_i},p_i(T^{r_i}))$.

\a \begin{boxedminipage}{160mm}
    \begin{algorithmic}[1]
		
	 \STATE\label{Alg1Step1}{\it STEP 1)} \ Apply the $\ast$-procedure to $f$ and go on with the $f^{(i_1,j_1)...(i_n,j_n)}$'s thus obtained. Form the graph $\textbf{G}_f$, applying the Stop Criterion.

	 \a \STATE\label{Alg1Step2}{\it STEP 2)} For any descending path $\g$ in $\textbf{G}_f$, work as in \ref{paths}, until you get a $p_{\g}(x) \in \C [x^{1\over r_\g}]$ which yields an approximation of a branch of $\bC$ at $O$.

 \a \STATE\label{Alg1Step3}{\it STEP 3)} \ Consider all $p_\g$'s obtained at the previous step and find the ones that give equivalent parameterizations. If $p_\gamma$ and $p_{\gamma '}$ are the truncations of Puiseux series $\tilde{p}_\gamma$ and $\tilde{p}_{\gamma '}$, they give equivalent parameterizations if and only if they have the same lowest common denominator $r$ and it exists $\omega_r\in \mathbb{C}$ with $\omega_r^r=1$, such that  $\tilde{p}_\gamma(x^{\frac{1}{r}}) = \tilde{p}_{\gamma '}(\omega_r x^{\frac{1}{r}})$ up to $\overline{n}$, where $\overline{n}=\min\left\lbrace n\in\N\; |\; f_n(x,y_n)\text{ and }f_n'(x,y_n)\text{ satisfy Theorem \ref{ift}}\right\rbrace $.

 Keep only one of the $p_\g$'s for each class of equivalence.

 \a \STATE \label{Alg1Step4}{\it STEP 4)} The number $s$ of equivalence classes found at {\it STEP 3)} is the number of branches of $\bC$ at $O$.  The data  $r_1,...,r_s$ and $p_{1} \in \C [x^{1\over r_{1}}],,...,p_s\in \C [x^{1\over r_{s}}]$ obtained before are the required output, and  $(T^{r_i},p_i(T^{r_i}))$, $i=1,...,s$ give approximations of the $s$ branches of our curve. $\textbf{END}$

    \end{algorithmic}
\end{boxedminipage}
\end{algorithm}

\medskip

\begin{algorithm}[H]\caption{Study of the branches at $(0,0)$ of a plane curve $\bC $ of equation $f=\sum_{i+j=1}^da_{ij}x^iy^j$, $a_{ij}\in \C$, $f=f_1^{n_1}...f_m^{n_m}$ (the irreducible decomposition of $f$), with $f_\ell \neq x$, $ \ell =1,...,m$}\label{Algorithm2}
\a  \textbf{Input}: $\bC $: $f=\sum_{i+j=1}^da_{ij}x^iy^j = f_1^{n_1}...f_m^{n_m}$, $O=(0,0)$.\\
 \textbf{Output}: Integers $r_{\ell, 1},...,r_{\ell, s_\ell}$ and $p_{\ell, 1} \in \C [x^{1\over r_{\ell, 1}}], \dots, p_{\ell, s_\ell} \in \C [x^{1\over r_{\ell, s_\ell}}]$, $\ell =1,...,m'\leq m$, such that $\bC $ possesses $s=s_1+...+ s_{m'}$ branches at $O$, and each branch is approximated by $(T^{r_{\ell, i}},p_{\ell,i}(T^{r_{\ell, i}}))$, $i\in \{1,...,{s_\ell}\}$.
\a \begin{boxedminipage}{160mm}
    \begin{algorithmic}[1]
		
	 \STATE\label{Alg2Step1}{\it STEP 1)} \ For each $f_\ell$, $\ell=1,...,m$, consider $f_\ell(O)$; let (up to reordering the $f_\ell$'s) $f_1(O)=...=f_{m'}(O)=0$,  while $f_\ell(O)\neq 0$ for $\ell >m'$. 
 
	 \a \STATE\label{Alg2Step2}{\it STEP 2)}  Apply Algorithm 1 to any $f_\ell$, $1 \leq \ell \leq m'$, so to get  $p_{\ell, 1} \in \C [x^{1\over r_{\ell, 1}}], \dots, p_{\ell, s_\ell} \in \C [x^{1\over r_{\ell, s_\ell}}]$  such that $(T^{r_{\ell, i}},p_{\ell,i}(T^{r_{\ell, i}}))$ are the (non-equivalent) approximations of the $s_\ell$ branches of $\bC_\ell : \{f_\ell =0\}$ at $O$.

 \a \STATE\label{Alg2Step3}{\it STEP 3)} \ Consider all $(T^{r_{\ell, i}},p_{\ell,i}(T^{r_{\ell, i}}))$, $i\in \{1,...,{s_\ell}\}$, $\ell \in \{1,...,m'\}$ obtained in the previous step. They give the approximations of all the branches of $\bC$ at $O$; each of them must be considered $n_\ell$ times.

 \a \STATE $\textbf{END}$

    \end{algorithmic}
\end{boxedminipage}

\end{algorithm}
 \subsection{The Lemmata}\label{lemmi}
 This subsection is of a more technical sort. We give the lemmata which justify why the algorithm actually works. 
\begin{lem}\label{equivalent}
	Let $\bC:f(x,y)=0$ be a reduced curve such that $x \nmid f$ and let 
	$$p(x^{1\over r})=\sum_{i=0}^\infty c_i(x^{1\over r})^i, \qquad p'(x^{1\over r'})=\sum_{i=0}^\infty c_i'(x^{1\over r'})^{i}$$
	be two roots of $f$. Moreover, let $\gamma$ and $\gamma'$ be the paths associated to $p$ and $p'$ with the respective notations  and set $\overline{n}=\min\left\lbrace n\in\N\; |\; f_n(x,y_n)\text{ and }f_n'(x,y_n)\text{ satisfy Theorem \ref{ift}}\right\rbrace $. Then the parameterizations
	$$(T^r,p(T)),\qquad (T^{r'},p'(T))$$
	are equivalent if and only if $r=r'$ and it exists $\omega$ such that $\omega^r=1$ and
	$$\sum_{i=0}^{\overline{n}}c_i(\omega x^{1\over r})^{i}=\sum_{i=0}^{\overline{n}}c_i'(x^{1\over r})^{i}.$$
\end{lem}
\proof  If the parameterizations are equivalent we conclude by Remark \ref{equiv}. By \cite{W}, Ch. IV, Theorem 4.1 it follows that since $p(x^{1\over r})$ is a root $f(x,y)$ then $p(\omega x^{1\over r})$ is a root of $f(x,y)$, where $\omega^r=1$. Moreover, since $f'_{\overline{n}}$ satisfies the Theorem \ref{ift}, it exists a unique root $q(x^{1\over r})$ of $f(x,y)$ such that $$q(x^{1\over r})=\sum_{i=0}^{\overline{n}}c_i'(x^{1\over r})^{i}+\dots $$
hence this root has to be $p'(x^{1\over r})$. 
By assumption, we have
$$p(\omega x^{1\over r})=\sum_{i=0}^\infty c_i(\omega x^{1\over r})^{i}=\sum_{i=0}^{\overline{n}}c_i(\omega x^{1\over r})^{i}+\dots=\sum_{i=0}^{\overline{n}}c_i'(x^{1\over r})^{i}+\dots$$
so that $p(\omega x^{1\over r})=q(x^{1\over r})$. Hence we have $p(\omega x^{1\over r})=p'(x^{1\over r})$ and the result follows by Remark \ref{equiv}.\prfend

\begin{lem}\label{lemma1} Notations as in subsection \ref{paths}. Let  {\rm $\textbf{a}_n\in\Gamma(f_n)$}, let $\rho$ denotes its slope, and let $c_1\dots c_t$ be the distinct roots of $g_n(1,y_n)$, of multiplicity $m_1,\dots,m_t$. Then $g_n(O)=0$, and the roots of $g_n(x,y_n)\in\C\{ \{x\}\}[y_n] $ are $q_{n1}(x)=c_1x^{r_n},\; \dots,\; q_{nt}(x)=c_tx^{r_n}$, where $r_n=-\frac{1}{\rho}$; more precisely,
 there is a constant $\alpha$ such that 
$$g_n(x,y_n)=\alpha \prod_{s=1}^{t}(y_n-c_sx^{r_n})^{m_s}$$
where $\sum_{i=1}^tm_i$ is equal to the height of {\rm$\textbf{a}_n$}. We say that  $m_s$ is {\em the multiplicity of $c_s$ as a root of $g_n$.}
\end{lem}

\proof We can assume that $y \nmid f_n$, hence $f_n$ is Newton convenient and we can consider its Newton Polygon $\Gamma(f_n)$; moreover, all its edges have strictly negative slopes by Remark \ref{novh}.
Hence, if $(i_1,j_1),\dots,(i_r,j_r)$ are the points of $\textbf{a}_n$ ordered from left to right, we have $i_1<i_2<\dots<i_r$ and $j_1>j_2>\dots>j_r$. So, if we set 
$$f_n(x,y_n)=\sum_{(i,j)\in\text{supp}(f_n)}a_{ij}x^iy_n^j$$ 
we have 
$$f_n^{\textbf{a}_n}=\sum_{k=1}^ra_{i_kj_k}x^{i_k}y_n^{j_k},\quad  u_i=i_1, v_i=j_s$$
which gives 
$$g_n(x,y_n)=\sum_{k=1}^{r}a_{i_kj_k}x^{i_k-i_1}y_n^{j_k-j_r}$$
so we see that $g_n(O)=0$. The slope of $\textbf{a}_n$ is $\rho=\frac{j_r-j_1}{i_r-i_1}$; we consider the possible roots of $g_n(x,y_n)\in\C\{ \{x\}\}[y_n] $ of the form $y_n=ux^{r_n}$, $u\in \C$, $r_n=-\frac{1}{\rho}$. We have
$$g_n(x,y_n)=g_n(x,ux^{r_n})=\sum_{k=1}^{r} a_{i_kj_k}x^{i_k-i_1+r_n(j_k-j_r)}u^{j_k-j_r}.$$ 
Since $(i_k,j_k)\in\textbf{a}_n,\;k=1,\dots,r$ the following equation holds
$$(j_1-j_r)(i_k-i_1)=(j_k-j_1)(i_1-i_r)\;\; k=1,\dots,r$$
and therefore 
$$
i_k-i_1+r_n(j_k-j_r)=\frac{(i_k-i_1)(j_r-j_1)-(j_k-j_r)(i_r-i_1)}{j_r-j_1}=i_r-i_1 \;\; k=1,\dots,r.
$$
Hence we have
$$
g_n(x,y_n)=x^{i_r-i_1}\sum_{k=1}^ra_{i_kj_k}u^{j_k-j_r}.
$$
Now we define $$h(y_n):=g_n(1,y_n)=\sum_{k=1}^{r}a_{i_kj_k}y_n^{j_k-j_r}.$$
Since $h(y_n)\in\C[y_n]$ and $\deg(h)=j_1-j_r$, it factorizes as 
$$h(y_n)=a_{i_1j_1}\prod_{s=1}^{t}(y_n-c_s)^{m_s},\qquad c_i\neq c_j\text{ and } \sum_{s=1}^{t}m_s=j_1-j_r.$$
So we have
$$
g_n(x,y_n)=x^{i_r-i_1}h(u)=a_{i_1j_1}x^{i_r-i_1}\prod_{s=1}^{t}(u-c_s)^{m_s}=a_{i_1j_1}x^{i_r-i_1}\prod_{s=1}^{t}(\frac{y_n}{x^{r_n}}-c_s)^{m_s}.
$$
Since $$r_n\sum_{s=1}^{t}m_s=r_n(j_1-j_r)=i_r-i_1$$
we can distribute the $x$ and we get
$$
g_n(x,y_n)=a_{i_1j_1}\prod_{s=1}^{t}x^{r_nm_s}(\frac{y_n}{x^{r_n}}-c_s)^{m_s}=a_{i_1j_1}\prod_{s=1}^{t}(y_n-c_sx^{r_n})^{m_s}.
$$
So the roots of $g_n$ are
$$q_{n1}=c_1x^{r_n},\; \dots,\; q_{nt}=c_tx^{r_n}$$
and we conclude the proof observing that $j_1-j_r$ is exactly the height of $\textbf{a}_n$.
\prfend

\begin{lem}\label{lemma2}  Notations as in subsection \ref{paths}. If  $cx^{r_n}$ is a root of $g_n(x,y_n)$ of multiplicity $m$, then the lowest pure power of $y_{n+1}$ in $f_{n+1} (x,y_{n+1})$ is $y_{n+1}^m$. Moreover, $f_{n+1} (O)=0$.
\end{lem}
\proof
We use the same notations introduced in the proof of Lemma \ref{lemma1}.  We can write $f_n$ as 
$$f_n=f_n^{\textbf{a}_n}+\underbrace{(f_n-f_n^{\textbf{a}_n})}_{\tilde{f_n}}=x^{i_1}y_n^{j_r}g_n+\tilde{f_n}.$$
Hence it follows that
$$f_n(x,x^{r_n}(c+y_{n+1}))=\underbrace{x^{i_1+r_n j_r}(c+y_{n+1})^{j_r}g_n(x,x^{r_n}(c+y_{n+1}))}_{A}+\tilde{f_n}(x,x^{r_n}(c+y_{n+1})).$$
By Lemma \ref{lemma1}
$$g_n(x,x^{r_n}(c+y_{n+1}))=a_{i_1j_1}\prod_{s=1}^{t}(x^{r_n} y_{n+1}+(c-c_s)x^{r_n})^{m_s}=a_{i_1j_1}x^{r_n\sum_{s=1}^{t}m_s}\prod_{s=1}^{t}(c-c_s+y_{n+1})^{m_s}$$
and
$$r_n=\frac{i_1-i_r}{j_r-j_1}\quad \text{and}\quad  \sum_{s=1}^{t}m_s=   j_1-j_r.$$
We can suppose without loss of generality that $c=c_1$ (and, as a consequence, that $m=m_1$). Therefore we have
$$g_n(x,x^{r_n}(c+y_{n+1}))=a_{i_1j_1}x^{ i_r-i_1}y_{n+1}^m\prod_{s=2}^{t}(c-c_s+y_{n+1})^{m_s}$$
and, by substituing this in $A$, we obtain
$$A=a_{i_1j_1}x^{ r_n j_r+i_r}y_{n+1}^m(c+y_{n+1})^{j_r}\prod_{s=2}^{t}(c-c_s+y_{n+1})^{m_s}.$$
Now let $a_{\overline{i}\overline{j}}x^{\overline{i}}y_n^{\overline{j}}$ be a monomial of $\tilde{f}_n$. Its evaluation in $(x,x^{r_n}(c+y_{n+1}))$ is 
$$a_{\overline{i}\overline{j}}x^{\overline{i}+r_n\overline{j}}(c+y_{n+1})^{\overline{j}}.$$ 
Since $(\bar i, \bar j)\in(\text{supp}(f_n)\setminus\textbf{a}_n)$, the point $(\bar i, \bar j)$ does not belong to $\textbf{a}_n$. Let  $\tilde a_n$ be the line containing the edge $\textbf{a}_n$; $\tilde a_n$ has equation $i+r_n j=i_r+r_n j_r$. Since $\Gamma (f_n)$ is part of the border of the convex hull of $\text{supp}(f_n)$, the point $(\bar i, \bar j)$ is in the half plane $ i+r_n j>i_r+r_n j_r$, hence
$$\bar i+r_n\bar j>i_r+r_n j_r. \eqno{(*)}$$ 
This shows that the power of $x$ for which we have to divide $f_n(x,x^{r_n}(c+y_{n+1}))$ is $x^{ r_n j_r+i_r}$ and that the pure powers of $y_{n+1}$ in $f_{n+1}$ are uniquely the ones coming from $$\frac{A}{x^{r_n j_r+i_r}}=a_{i_1j_1}y^m_{n+1}(c+y_{n+1})^{j_r}\prod_{s=2}^{t}(c-c_s+y_{n+1})^{m_s}.$$
Hence the monomial with the lowest pure power of $y_{n+1}$ is 
$$
\left( a_{i_1j_1}c^{j_r}\prod_{s=2}^{t}(c-c_s)^{m_s}\right) y_{n+1}^m
$$
so the lowest pure power of $y_{n+1}$ in $f_{n+1} (x,y_{n+1})$ is $y_{n+1}^m$.
We have $$f_{n+1}(x,y_{n+1})=\frac{A}{x^{r_n j_r+i_r}}+{\tilde{f_n}(x,x^{r_n}(c+y_{n+1}))\over x^{r_n j_r+i_r}}$$
and $\frac{A}{x^{r_n j_r+i_r}}(O)=0$ since it is divided by $y_{n+1}$, while ${\tilde{f_n}(x,x^{r_n}(c+y_{n+1}))\over x^{r_n j_r+i_r}}(O)=0$ by $(*)$.
\prfend
\begin{lem}\label{lemma3} Notations as in subsection \ref{paths}. Let $n\in \N$ be the minimum such that there exists $t \in \N$, $t \geq 1$ such that  $y_{n+1}^t\mid f_{n+1}$. Then: 
$$y_{n+1}^t\mid f_{n+1}\text{ in }\C\left\lbrace \left\lbrace x\right\rbrace \right\rbrace[y_{n+1}] \Leftrightarrow  (y-\sum_{i=0}^{n}c_ix^{r_i+r_{i-1}+\dots+r_0})^t\mid f\text{ in }\C\left\lbrace  \{x\} \right\rbrace [y]$$
\end{lem}
\proof Let
$$
f_n(x,y_n)=\sum_{(i,j)\in\text{supp}(f_n)}a_{ij}x^iy_n^j\;\;;\qquad f_{n+1}(x,y_{n+1})=\sum_{(i',j')\in\text{supp}(f_{n+1})}b_{i'j'}x^{i'}y_{n+1}^{j'}
$$
and recall that the exponents of $x$ are rationals, while those of $y$ are integers.
\a $\Rightarrow$) Since $y_{n+1}^t\mid f_{n+1}(x,y_{n+1})$ it exists $h_{n+1}\in\C\{\{x\}\}[y_{n+1}]$ such that $f_{n+1}(x,y_{n+1})=y_{n+1}^th_{n+1}(x,y_{n+1})$. Moreover we have
$$f_n(x,y_n)=x^{m_n}f_{n+1}\left( x,\frac{y_n}{x^{r_n}}-c_n\right) =x^{m_n}\left( \frac{y_n}{x^{r_n}}-c_n\right) ^th_{n+1}\left( x,\frac{y_n}{x^{r_n}}-c_n\right)=$$
$$=(y_n-c_nx^{r_n})^t\underbrace{ x^{m_n-tr_n}h_{n+1}\left( x,\frac{y_n}{x^{r_n}}-c_n\right)}_{h_n(x,y_n)}.$$
Now we want to show that $h_n(x,y)\in\C\{\{x\}\}[y_n]$. We have
$$h_{n+1}(x,y_{n+1})=\frac{f_{n+1}(x,y_{n+1})}{y_{n+1}^t}=\sum_{(i',j')\in\text{supp}(f_{n+1})}b_{i'j'}x^{i'}y_{n+1}^{j'-t}$$
therefore
$$
h_n(x,y_n)=x^{m_n-tr_n}h_{n+1}\left( x,\frac{y_n}{x^{r_n}}-c_n\right) =x^{m_n-tr_n}\sum_{(i',j')\in\text{supp}(f_{n+1})}b_{i'j'}x^{i'}\left( \frac{y_n}{x^{r_n}}-c_n\right) ^{j'-t}=
$$
$$
=\sum_{(i',j')\in\text{supp}(f_{n+1})}b_{i'j'}x^{m_n-tr_n+i'-r_n(j'-t)}\left({y_n}-c_nx^{r_n}\right) ^{j'-t}.
$$
Since $j'-t\geq 0$, it is enough to prove that
$$m_n-tr_n+i'-r_n(j'-t)\geq0,\quad \forall\;\;(i',j')\in\text{supp}(f_{n+1}).\eqno{(\bullet)}$$
By definition
$$f_{n+1}(x,y_{n+1})=\frac{f_n(x,x^{r_n}(c_n+y_{n+1}))}{x^{m_n}}=\sum_{(i,j)\in\text{supp}(f_n)}a_{ij}x^{i+jr_n-m_n}(c_n+y_{n+1})^j=$$
$$=\sum_{(i,j)\in\text{supp}(f_n)}a_{ij}x^{i+jr_n-m_n}\sum_{k=0}^{j}{j \choose k}y_{n+1}^kc_n^{j-k}$$
hence
$$(i',j')\in\text{supp}(f_{n+1})\Rightarrow\exists\;\;(i,j)\in\text{supp}(f_n),k\in\N, k\leq j\; \text{such that } i'=i+jr_n-m_n\text{ and } j'=k. $$ 
So, in order that $(\bullet)$ holds it suffices that
$$
m_n-tr_n+i+jr_n-m_n-r_n(k-t)=i+(j-k)r_n\geq0\quad \forall\;\; (i,j)\in\text{supp}(f_n), 0\leq k\leq j
$$
and this inequality holds because $i,(j-k),r_n\geq0$. Thus we showed that $h_n(x,y_n)\in\C\{\{x\}\}[y]$.
\dd Now, by the same reasoning, we have
$$f_{n-1}(x,y_{n-1})=x^{m_{n-1}}f_n\left( x,\frac{y_{n-1}}{x^{r_{n-1}}}-c_{n-1}\right)=x^{m_{n-1}}\left(\frac{y_{n-1}}{x^{r_{n-1}}}-c_{n-1}-c_{n}x^{r_n}\right) ^t h_n\left( x,\frac{y_{n-1}}{x^{r_{n-1}}}-c_{n-1}\right) = $$
$$
=(y_{n-1}-c_{n-1}x^{r_{n-1}}-c_n^{r_n+r_{n-1}})^t\underbrace{ x^{m_{n-1}-tr_{n-1}}h_n\left( x,\frac{y_{n-1}}{x^{r_{n-1}}}-c_{n-1}\right) }_{h_{n-1}(x,y_{n-1})\in\C\left\lbrace \left\lbrace x\right\rbrace \right\rbrace  [y_{n-1}]}.
$$
Hence, since $f=f_0$, iterating the same procedure we find $h(x,y)\in\C\left\lbrace \left\lbrace x\right\rbrace \right\rbrace [y]$ such that  
$$f(x,y)=\left( y-\sum_{i=0}^{n}c_ix^{r_i+r_{i-1}+\dots+r_0}\right) ^th(x,y).$$
$\Leftarrow$) The proof is analogous to the previous one.
\prfend
\begin{lem}\label{lemma35} 	Let us consider $f,g\in\C[u,v]$ such that $g$ is irreducible, $\deg g>0$, and let $t\in \N, t>0$. Then $g^t\mid f$ if and only if $g$ divides $f$ and all its partial derivatives up to order $t-1$.	
\end{lem}
\proof $\Rightarrow$) Is trivial.\\
$\Leftarrow$)  Let us prove the statement for $t=1$. Let $f=gh$; since $g\mid f_u$ and $g\mid f_v$, there exist $h_1,h_2\in\C[u,v]$ such that $f_u=gh_1$ and $f_v=gh_2$. We have
$$gh_1=f_u=(gh)_u=g_uh+gh_u, \qquad gh_2=f_v=(gh)_v=g_vh+gh_v.$$
Thus, we have that $g(h_1-h_u)=g_uh$ and $g(h_2-h_v)=g_vh$ so that $g\mid g_uh$ and $g\mid g_vh$. Since $g$ is irreducible it follows that either $g\mid g_u$ or $g\mid h$ and either $g\mid g_v$ or $g\mid h$. Since $\deg g_u<\deg g$, we have that $g\mid g_u$ if and only if $g_u=0$, and similarly for $g_v$. Now, notice that since $\deg g>0$, at least one between $g_u$ and $g_v$ is not the zero polynomial and hence $g\mid h$. \\
Now let us suppose that  $t>1$ and that $g$ divides $f$ and all its derivatives up to order $t-1$. In particular $g$ divides all derivatives of order $t-2$ and their first partial derivatives. Thus, by the previous part of the proof, it follows that $g^2$ divides all partial derivatives of $f$ of order $t-2$. By iterating this argument we have that $g^a$ divides all partial derivatives of $f$ of order $t-a$ and hence $g^t\mid f$. 
\prfend
\begin{lem}\label{lemma4}
	 Notations as in subsection \ref{paths}. Let $t, n \in \N$, $t \geq 1, n\geq 0$, and let \linebreak $y_{n+1}^t\mid f_{n+1}$. Then if $r$ is the lowest common denominator of $\sum_{i=0}^{n}c_ix^{r_i+r_{i-1}+\dots+r_0}$ we have that $\bR=(T^r,\sum_{i=0}^{n}c_iT^{r\left( r_i+r_{i-1}+\dots+r_0\right) })$ is a branch of $\bC:f=0$ at $O$.  Moreover, if  $g(x,y)=0$ is an equation of the parametric curve $\bR$, then $g^t\mid f$ in $\C[x,y]$. \end{lem}
\proof 
We use the same notations introduced in the proof of Lemma \ref{lemma3}.
\b By assumption $y_{n+1}^t\mid f_{n+1}$, hence by Lemma \ref{lemma3} we have  that there exists $h \in \C\left\lbrace  \{x\} \right\rbrace [y]$ such that
$f=h\cdot (y-\sum_{i=0}^{n}c_ix^{r_i+r_{i-1}+\dots+r_0})^t$. In particular $$f(x,\sum_{i=0}^{n}c_ix^{r_i+r_{i-1}+\dots+r_0})=0$$ and hence $\bR$ is a branch of $\bC$ at $O$. 
\dd Since $\bR$ has a polynomial parameterization, it is an algebraic rational curve and therefore there exists $g\in\C[x,y]$, $g$ irreducible, such that $\bR:g=0$. In order to prove that $g^t\mid f$ in $\C[x,y]$, we start proving that $g^t\mid f$ in $\C\{\{x\}\}[y]$.  Notice that, by the proof of Lemma \ref{lemma3}, it follows that there exists $m\in\N$ such that $h, (y-\sum_{i=0}^{n}c_ix^{r_i+r_{i-1}+\dots+r_0})\in\C[x^{1\over m},y]$. As a consequence we can prove that $g^t\mid f$ in $\C[x^{1\over m},y]$. To do this, we apply Lemma \ref{lemma35} to $f,g\in\C[x^{1\over m},y]$. Let  $$P=(\overline{T}^r,\sum_{i=0}^{n}c_i\overline{T}^{r(r_i+r_{i-1}+\dots+r_0)})\in\bR,$$ then 
$$f(P)=f(\overline{T}^r,\sum_{i=0}^{n}c_i\overline{T}^{r(r_i+r_{i-1}+\dots+r_0)})=0,\text{ }\forall P\in\bR\;\Rightarrow\; P\in\bC\text{ }\forall P\in\bR\;\Rightarrow\; g\mid f. $$
In order to calculate the partial derivatives of $f$ with respect to $x^{1\over m}$ and $y$, that we denote respectively by $f_{x^{1\over m}}$ and $f_y$ , we set $a_i:=mr_i \in \Z$ and rewrite $f$ as follows
$$f(x^{1\over m},y)=(y-\sum_{i=0}^{n}c_i(x^{1\over m})^{a_i+a_{i-1}+\dots+a_0})^t\cdot h(x^{1\over m},y)$$
so that the first order partial derivates of $f$ are
$$f_{x^{1\over m}}=-t(\sum_{i=0}^{n}c_i(a_i+a_{i-1}+\dots+a_0)(x^{1\over m})^{a_i+a_{i-1}+\dots+a_0-1})(y-\sum_{i=0}^{n}c_i(x^{1\over m})^{a_i+a_{i-1}+\dots+a_0})^{t-1}h(x^{1\over m},y)+$$
$$+(y-\sum_{i=0}^{n}c_i(x^{1\over m})^{a_i+a_{i-1}+\dots+a_0})^th_{x^{1\over m}}(x^{1\over m},y) $$
and
$$f_y=t(y-\sum_{i=0}^{n}c_i(x^{1\over m})^{a_i+a_{i-1}+\dots+a_0})^{t-1}h(x^{1\over m},y)+(y-\sum_{i=0}^{n}c_i(x^{1\over m})^{a_i+a_{i-1}+\dots+a_0})^th_y(x^{1\over m},y).$$
So, we have
$$f_{x^{1\over m}}(P)=f_y(P)=0\;\;\forall P\in\bR\Rightarrow g\mid f_{x^{1\over m}}\text{ and }g\mid f_y.$$
In the same way, each partial derivative up to $t-1$ order is sum of products that have at least one factor $(y-\sum_{i=0}^{n}c_i(x^{1\over m})^{a_i+a_{i-1}+\dots+a_0})$, thus all of them vanish in each $P\in\bR$. As a consequence, by Lemma \ref{lemma35} we have that $g^t\mid f$ in $\C[x^{1\over m},y]$, that is there exists $p\in\C[x^{1\over m},y]$ such that $f=g^tp$. Since $f,g\in\C[x,y]\subseteq\C[x^{1\over m},y]$, in order to show that $g^t\mid f$ in $\C[x,y]$ it is enough to show that $p\in\C[x,y]$. Let us write $p$ as $p=p_1+p_2$, where each monomial of $p_1$ is in $\C[x,y]$ and each monomial of $p_2$ is in $\C[x^{1\over m},y]\setminus\C[x,y]$.  We have that $$g^tp_2=f-g^tp_1$$
so that, since $f,g,p_1\in\C[x,y]$, also $g^tp_2\in\C[x,y]$. Now we observe that if $x^ay^b$ is a monomial of $g^tp_2$ then $a$ is sum of an integer, coming from a monomial of $g^t$, and of a rational that is no an integer, coming from a monomial of $p_2$, and thus $a$ is a rational but not an integer. As a consequence we have that $g^tp_2=0$. Hence, since $\C[x^{1\over m},y]$ is a domain and $g\neq 0$, we have that $p_2=0$ and $p=p_1\in\C[x,y]$ and this ends the proof.\prfend

\subsection{An example} 
The following example illustrates how the Algorithm works. 
\begin{exa}\rm 
	Let consider the integral curve $\bC:f=0$, where
	$$f=2y^6+6xy^5-8x^3y^3+2x^3y^4+(2\sqrt{3}+2)x^4y^3+(4\sqrt{3}-4)x^5y^2+(\sqrt{3}-2)x^7y+\frac{\sqrt{3}-2}{8}x^{10}+2x^{11}.$$
We have $m_O(\bC)=6$. The Newton polygon $\Gamma(f_0)=\Gamma(f)$ is the following:
	\begin{center}
		\includegraphics{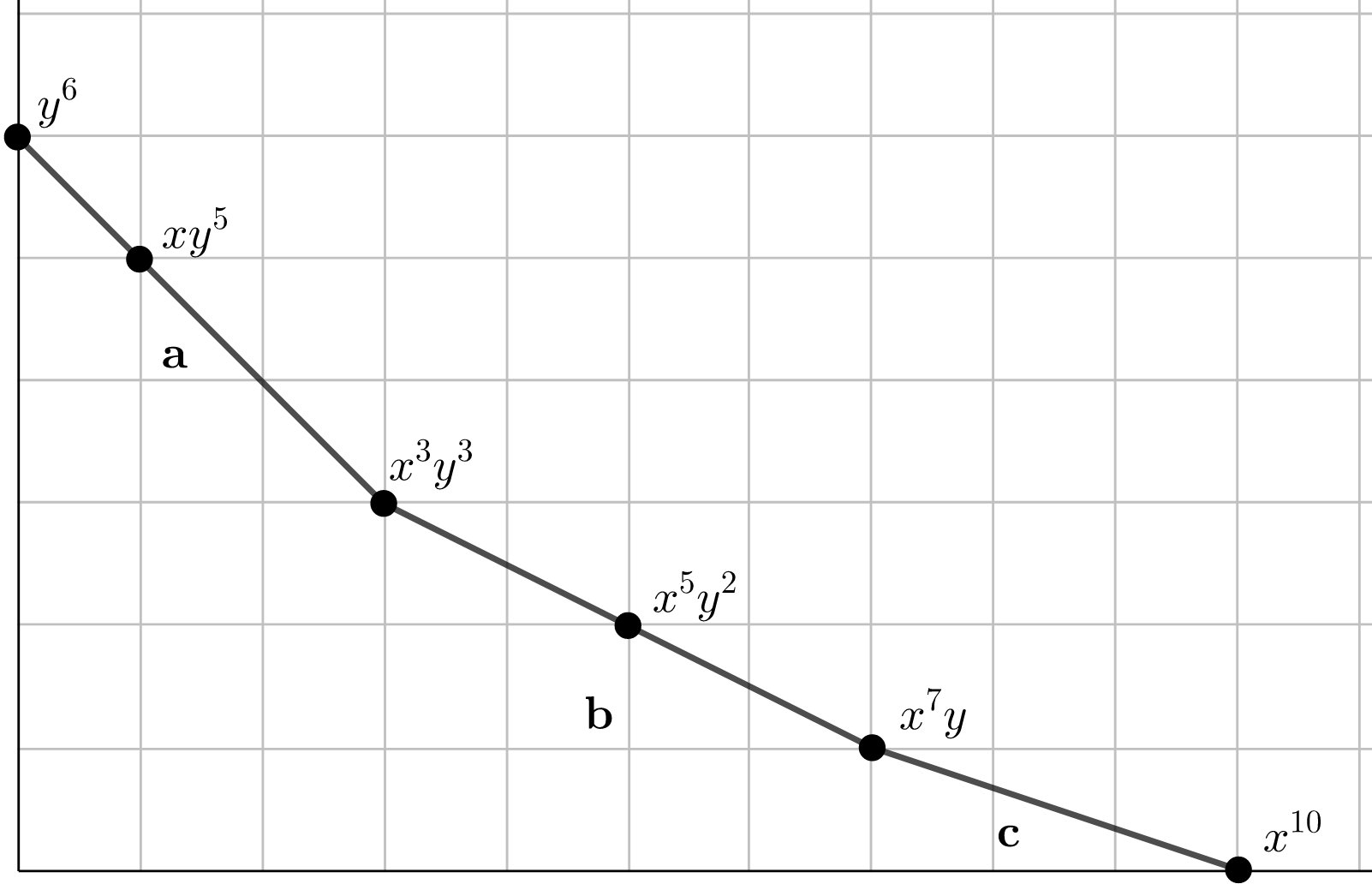}
	\end{center}
 We start setting $f_0=f$. Let us analyze each possible choice individually:
	\begin{itemize}[leftmargin=*]
		\item $\textbf{a}_0=\textbf{a}$\\
		We have
		$$f_0^{\textbf{a}_0}=f^{\textbf{a}_0}=2y^6+6xy^5-8x^3y^3$$
		$$g_0=\frac{f_0^{\textbf{a}_0}}{y^3}=2(y^3+3xy^2-4x^3)=2(y+2x)^2(y-x).$$
		So we have two distinct roots of $g_0$, namely $q_{01}=-2x$ double and $q_{02}=x$ simple, and, as a consequence, two possible choice for $q_0$:
		\renewcommand\labelitemii{ \textasteriskcentered }
		\begin{itemize}[leftmargin=*]
			\item $q_0=q_{01}$\\
			We have 
			$$f_0(x,x(-2+y_1))=x^6(48y_1^2-88y_1^3+60y_1^4-18y_1^5+2y_1^6+(8\sqrt{3}-24)xy_1+(-8\sqrt{3}+32)xy_1^2+$$
			$$+(2\sqrt{3}-14)xy_1^3+2xy_1^4+(-2\sqrt{3}+4)x^2+(\sqrt{3}-2)x^2y_1+(\frac{\sqrt{3}}{8}-\frac{1}{4})x^4+2x^5 ) $$
			hence
			$$f_1(x,y_1)=\frac{f_0(x,x(-2+y_1))}{x^6}=48y_1^2-88y_1^3+60y_1^4-18y_1^5+2y_1^6+(8\sqrt{3}-24)xy_1+(-8\sqrt{3}+32)xy_1^2+$$
			$$+(2\sqrt{3}-14)xy_1^3+2xy_1^4+(-2\sqrt{3}+4)x^2+(\sqrt{3}-2)x^2y_1+(\frac{\sqrt{3}}{8}-\frac{1}{4})x^4+2x^5.$$
			The Newton polygon $\Gamma(f_1)$ is the following
			\begin{center}
				\includegraphics{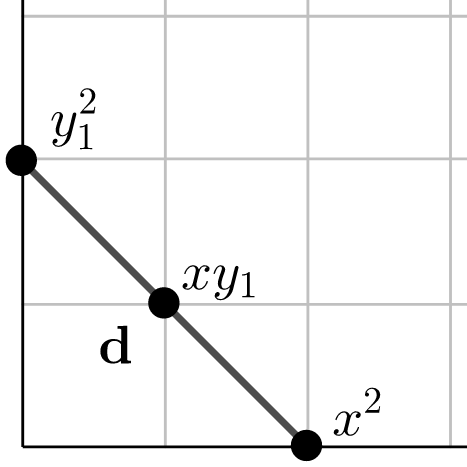}
			\end{center}
			We have an unique edge $\textbf{a}_1=\textbf{d}$, therefore we get
			$$g_1=f_1^{\textbf{a}_1}=48y_1^2+(8\sqrt{3}-24)xy_1+(-2\sqrt{3}+4)x^2=(4\sqrt{3}y+(-\sqrt{3}+1)x)^2.$$
			There is an unique root of $g_1$, which is
			$$q_1=-\frac{-\sqrt{3}+1}{4\sqrt{3}}x=\frac{3-\sqrt{3}}{12}x.$$
			So we have   
			$$f_1(x,x(\frac{3-\sqrt{3}}{12}+y_2))=x^2(48y_2^2+(\frac{-\sqrt{3}}{18}+\frac{1}{6})x+(\frac{8\sqrt{3}}{3}-4)xy_2+(-\frac{\sqrt{3}}{72}-\frac{1}{16})x^2+(14\sqrt{3}-34)xy_2^2+$$
			$$+(-\frac{23\sqrt{3}}{12}+\frac{13}{4})x^2y_2+(\frac{41\sqrt{3}}{1152}+\frac{2233}{1152})x^3-88xy_2^3+(-10\sqrt{3}+18)x^2y_2^2+(\frac{35\sqrt{3}}{72}-\frac{27}{32})x^3y_2+$$
			$$+(-\frac{5\sqrt{3}}{2304}+\frac{13}{3456})x^4+(-18\sqrt{3}+46)x^2y_2^3+(\frac{21\sqrt{3}}{8}-\frac{37}{8})x^3y_2^2+(-\frac{19\sqrt{3}}{576}+\frac{11}{192})x^4y_2+60x^2y_2^4+$$
			$$+(\frac{41\sqrt{3}}{6}-13)x^3y_2^3+(-\frac{5\sqrt{3}}{24}+\frac{35}{96})x^4y_2^2+(\frac{15\sqrt{3}}{2}-\frac{41}{2})x^3y_2^4+(-\frac{25\sqrt{3}}{36}+\frac{5}{4})x^4y_2^3-18x^3y_2^5+$$
			$$+(-\frac{5\sqrt{3}}{4}+\frac{5}{2})x^4y_2^4+(-\sqrt{3}+3)x^4y_2^5+2x^4y_2^6)$$
			$$f_2(x,y_2)=\frac{f_1(x,x(\frac{3-\sqrt{3}}{12}+y_2))}{x^2}.$$
			The Newton polygon $\Gamma(f_2)$ is the following
			\begin{center}
				\includegraphics{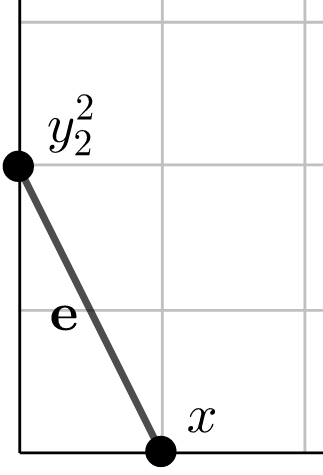}
			\end{center}We have an unique edge $\textbf{a}_2=\textbf{e}$ and we get:
			$$g_2=f_2^{\textbf{a}_2}=48y_2^2+(\frac{-\sqrt{3}}{18}+\frac{1}{6}x)$$
			whose roots are
			$$q_{21}=i\sqrt{\frac{3-\sqrt{3}}{864}}x^{\frac{1}{2}}\qquad q_{22}=-i\sqrt{\frac{3-\sqrt{3}}{864}}x^{\frac{1}{2}}.$$
			By Lemma \ref{lemma2}, choosing either $q_2=q_{21}$ or $q_2=q_{22}$, the polynomial $f_3$ satisfies the hypothesis of Theorem \ref{ift}; hence by \ref{stop2} we can stop. Since $$r_0=1,\quad r_1=1, \quad r_2={1\over 2}$$ we get
			$$ p_1=-2x+\frac{3-\sqrt{3}}{12}x^2+i\sqrt{\frac{3-\sqrt{3}}{864}}x^{\frac{5}{2}}+\dots$$
			$$ p_2=-2x+\frac{3-\sqrt{3}}{12}x^2-i\sqrt{\frac{3-\sqrt{3}}{864}}x^{\frac{5}{2}}+\dots.$$
			
			Therefore we have the following parameterizations
			$$R_1=(T^2,-2T^2+\frac{3-\sqrt{3}}{12}T^4+i\sqrt{\frac{3-\sqrt{3}}{864}}T^5+\dots)$$ $$R_2=(T^2,-2T^2+\frac{3-\sqrt{3}}{12}T^4-i\sqrt{\frac{3-\sqrt{3}}{864}}T^5+\dots).$$
			We have that $p_1(x) = p_2(-x)$, at $O$, so by Step 3 of the Algorithm with $\omega_2 = -1$ we conclude that $R_1$ and $R_2$ are equivalent parameterizations of the same 2-branch of $\bC$.
			\item $q_0=q_{02}$\\
			We have
			$$f_0(x,x(1+y_1))=x^6(18y_1+6\sqrt{3}x+66y_1^2+(14\sqrt{3}+6)xy_1+(\sqrt{3}-2)x^2+92y_1^3+(10\sqrt{3}+14)xy_1^2+$$
			$$+(\sqrt{3}-2)x^2y_1+60y_1^4+(2\sqrt{3}+10)xy_1^3+(\frac{\sqrt{3}}{8}-\frac{1}{4})x^4+18y_1^5+2xy_1^4+2x^5+2y_1^6)$$
			and hence
			$$f_1(x,y_1)=\frac{f_0(x,x(1+y_1))}{x^6}=18y_1+6\sqrt{3}x+66y_1^2+(14\sqrt{3}+6)xy_1+(\sqrt{3}-2)x^2+92y_1^3+$$
			$$+(10\sqrt{3}+14)xy_1^2+(\sqrt{3}-2)x^2y_1+60y_1^4+(2\sqrt{3}+10)xy_1^3+(\frac{\sqrt{3}}{8}-\frac{1}{4})x^4+18y_1^5+2xy_1^4+2x^5+2y_1^6).$$
			The Newton polygon $\Gamma(f_1)$ is the following
			\begin{center}
				\includegraphics{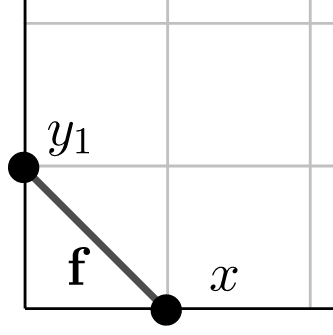}
			\end{center}
			We can just choose $\textbf{a}_1=\textbf{f}$ getting
			$$g_1=f_1^{\textbf{f}}=18y_1+6\sqrt{3}x$$
			whose unique root is
			$$q_1=-\frac{\sqrt{3}}{3}x.$$
			Moreover $f_1$ satysfies the hypothesis of Theorem \ref{ift}, therefore by Remark \ref{stop2} we can stop by getting the following root of $f$
			$$p_3(x)=x-\frac{\sqrt{3}}{3}x^2+\dots.$$
			Hence we have the parameterization
			$$R_3=(T,T-\frac{\sqrt{3}}{3}T^2+\dots)$$
			that is a 1-branch of $\bC$ at $O$.
		\end{itemize}
		\item $\textbf{a}_0=\textbf{b}$\\
		We have
		$$f_0^{\textbf{a}_0}=f^{\textbf{a}_0}=-8x^3y^3+(4\sqrt{3}-4)x^5y^2+(\sqrt{3}-2)x^7y$$
		$$g_0=\frac{f_0^{\textbf{a}_0}}{x^3y}=-8y^2+(4\sqrt{3}-4)x^2y^2+(\sqrt{3}-2)x^4=-(2\sqrt{2}y+\frac{\sqrt{2}-\sqrt{6}}{2}x^2)^2.$$
		Hence we have only one choice for $q_0$ that is
		$$q_0=-\frac{\sqrt{2}-\sqrt{6}}{4\sqrt{2}}x^2=\frac{\sqrt{3}-1}{4}x^2.$$
		Proceeding we get
		$$f_0(x,x^2(\frac{\sqrt{3}-1}{4}+y_1))=x^9((-2\sqrt{3}+2)y_1^2+(\frac{3\sqrt{3}}{4}-\frac{3}{4})xy_1+(\frac{17\sqrt{3}}{128}+\frac{227}{128})x^2-8y_1^3+3xy_1^2+$$
		$$+(-\frac{9\sqrt{3}}{8}+\frac{65}{32})x^2y_1+(\frac{-15\sqrt{3}}{256}+\frac{13}{128})x^3+(2\sqrt{3}+2)xy_1^3+(\frac{33\sqrt{3}}{8}-\frac{51}{8})x^2y_1^2+(\frac{33\sqrt{3}}{64}-\frac{57}{64})x^3y_1+$$
		$$+(-\frac{11\sqrt{3}}{2}+13)x^2y_1^3+(-\frac{15\sqrt{3}}{8}+\frac{105}{32})x^3y_1^2+(\frac{15\sqrt{3}}{2}-\frac{11}{2})x^2y_1^4+(\frac{15\sqrt{3}}{4}-\frac{25}{4})x^3y_1^3+6x^2y_1^5+$$
		$$+(-\frac{15\sqrt{3}}{4}+\frac{15}{2})x^3y_1^4+(3\sqrt{3}-3)x^3y_1^5+2x^3y_1^6)$$
		$$f_1(x,y_1)=\frac{f_0(x,x^2(\frac{\sqrt{3}-1}{4}+y_1))}{x^9}.$$
		The Newton Polygon $\Gamma(f_1)$ is the following
		\begin{center}
			\includegraphics{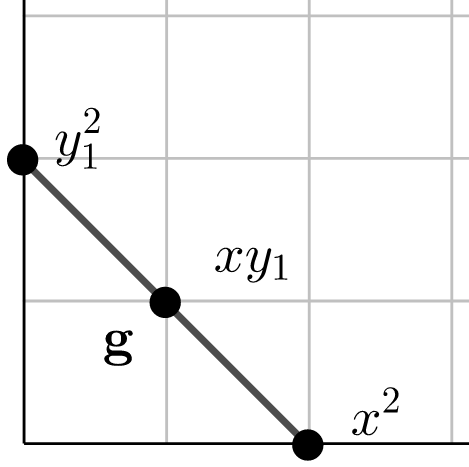}
		\end{center}
		We can only choose $\textbf{a}_1=\textbf{g}$ getting
		$$g_1=f_1^{\textbf{g}}=(-2\sqrt{3}+2)y_1^2+(\frac{3\sqrt{3}}{4}-\frac{3}{4})xy_1+(\frac{17\sqrt{3}}{128}+\frac{227}{128})x^2$$
		whose roots are
		$$q_{11}=\frac{3-3\sqrt{3}+2\sqrt{35+48\sqrt{3}}}{16(-\sqrt{3}+1)}x$$
		$$q_{12}=\frac{3-3\sqrt{3}-2\sqrt{35+48\sqrt{3}}}{16(-\sqrt{3}+1)}x.$$
		By Lemma \ref{lemma2}, choosing either $q_1=q_{11}$ or $q_1=q_{12}$, the polynomial $f_2$ satisfies the hypothesis of Theorem \ref{ift}. Hence, by Remark $\ref{stop2}$ we can stop getting the following roots of $f$
		$$p_4=\frac{\sqrt{3}-1}{4}x^2+\frac{3-3\sqrt{3}+2\sqrt{35+48\sqrt{3}}}{16(-\sqrt{3}+1)}x^3+\dots$$
		$$p_5=\frac{\sqrt{3}-1}{4}x^2+\frac{3-3\sqrt{3}-2\sqrt{35+48\sqrt{3}}}{16(-\sqrt{3}+1)}x^3+\dots.$$
		Hence we have the parameterizations
		$$R_4=(T,\frac{\sqrt{3}-1}{4}T^2+\frac{3-3\sqrt{3}+2\sqrt{35+48\sqrt{3}}}{16(-\sqrt{3}+1)}T^3+\dots)$$
		$$R_5=(T,\frac{\sqrt{3}-1}{4}T^2+\frac{3-3\sqrt{3}-2\sqrt{35+48\sqrt{3}}}{16(-\sqrt{3}+1)}T^3+\dots)$$
		that give two distinct 1-branches of	$\bC$ at $O$.
		\item $\textbf{a}_0=\textbf{c}$\\
		We have 
		$$f_0^{\textbf{c}}=f^\textbf{c}=(\sqrt{3}-2)x^7y+\frac{\sqrt{3}-2}{8}x^{10}$$
		$$g_0=\frac{f_0^{\textbf{c}}}{x^7}=(\sqrt{3}-2)y+\frac{\sqrt{3}-2}{8}x^{3}.$$
		We have just a choice for $q_0$ that is
		$$q_0=-\frac{1}{8}x^3.$$
		Proceeding we get
		$$f_0(x,x^3(-\frac{1}{8}+y_1))=x^{10}((\sqrt{3}-2)y_1+(\frac{\sqrt{3}}{16}+\frac{31}{16}))x+(-\sqrt{3}+1)xy_1+\frac{1}{64}x^2+(4\sqrt{3}-4)xy_1^2-\frac{3}{8}x^2y_1+$$
		$$+(-\frac{\sqrt{3}}{256}-\frac{1}{256})x^3+3x^2y_1^2+(\frac{3\sqrt{3}}{32}+\frac{3}{32})x^3y_1-8x^2y_1^3+(-\frac{3\sqrt{3}}{4}-\frac{3}{4})x^3y_1^2+\frac{1}{2048}x^5+(2\sqrt{3}+2)x^3y_1^3+$$
		$$-\frac{1}{64}x^5y_1-\frac{3}{16384}x^6+\frac{3}{16}x^5y_1^2+\frac{15}{2048}x^6y_1-x^5y_1^3-\frac{15}{128}x^6y_1^2+\frac{1}{131072}x^8+2x^5y_1^4+\frac{15}{16}x^6y_1^3-\frac{3}{8192}x^8y_1+$$
		$$-\frac{15}{4}x^6y_1^4+\frac{15}{2048}x^8y_1^2+6x^6y_1^5-\frac{5}{64}x^8y_1^3+\frac{15}{32}x^8y_1^4-\frac{3}{2}x^8y_1^5+2x^8y_1^6$$
		$$f_1(x,y_1)=\frac{f_0(x,x^3(-\frac{1}{8}+y_1))}{x^{10}}.$$
		The Newton polygon $\Gamma(f_1)$ is the following
		\begin{center}
			\includegraphics{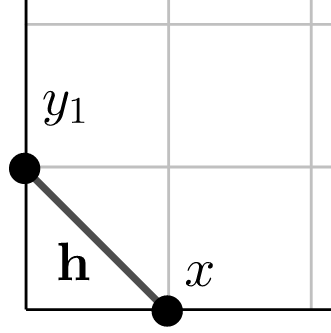}
		\end{center}
		We can only choose $\textbf{a}_1=\textbf{h}$ getting
		$$g_1=f_1^{\textbf{h}}=(\sqrt{3}-2)y_1+(\frac{\sqrt{3}}{16}+\frac{31}{16})x$$
		whose unique root is
		$$q_1=-\frac{33\sqrt{3}+65}{16}x.$$
		Moreover $f_1$ satysfies the hypothesis of Theorem \ref{ift}, therefore by Remark \ref{stop2} we can stop getting the following root of $f$
		$$p_6=-\frac{1}{8}x^3-\frac{33\sqrt{3}+65}{16}x^4+\dots$$
		which gives the parameterization
		$$R_6(T,-\frac{1}{8}T^3-\frac{33\sqrt{3}+65}{16}T^4+\dots)$$
		that is a 1-branch of $\bC$ at $O$. \\
		
	\end{itemize}
 We have studied all possible choices so we can stop the analysis.
\dd Since we saw that $R_1$ and $R_2$ are equivalent parameterizations of the same 2-branch of $\bC$, we have found five branches of $\bC$, as the following figure shows.
	\begin{center}
		\includegraphics[scale=2.2]{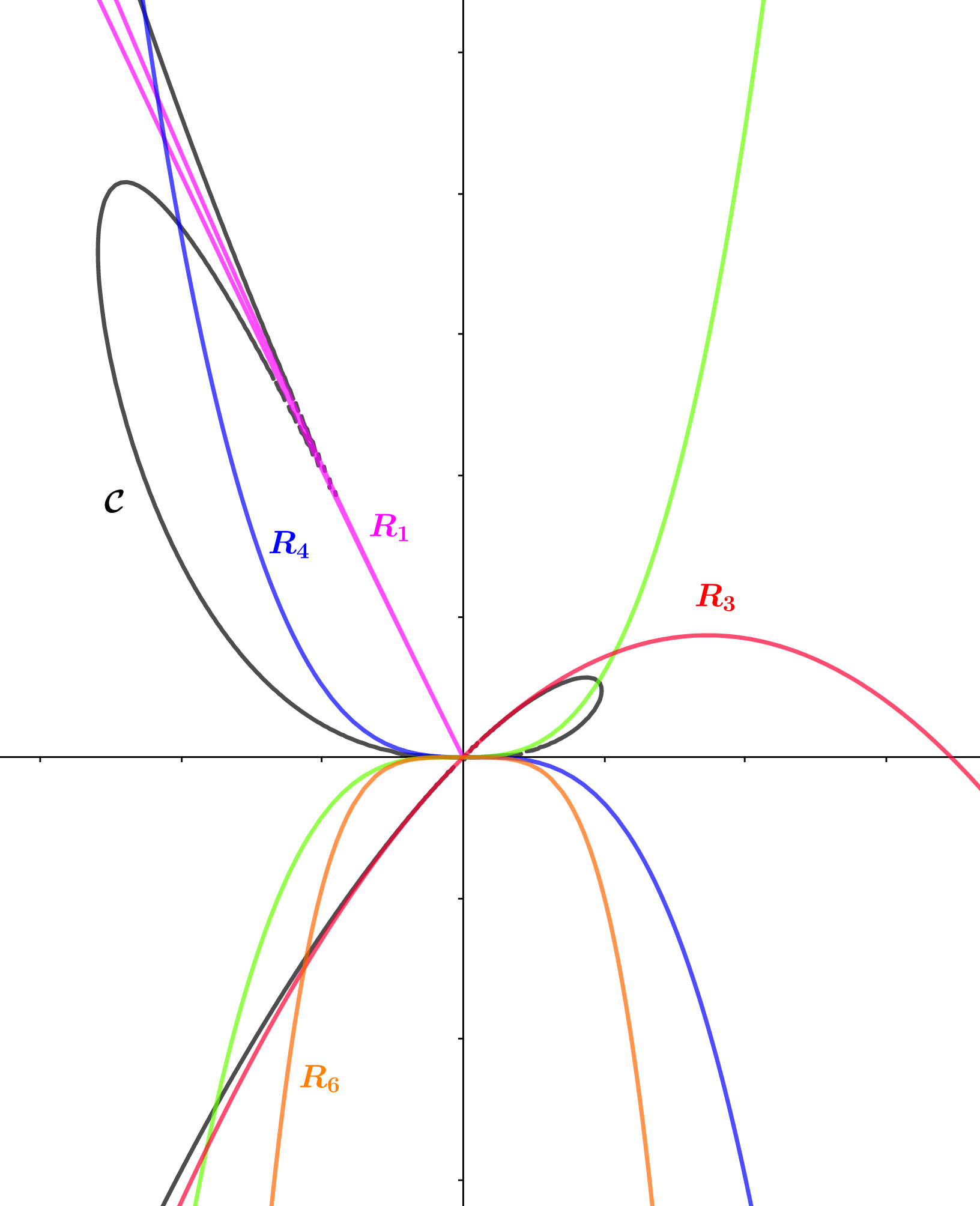}
	\end{center}
\end{exa}

\section{Triple points }\label{triple}
\subsection{Output of the Newton-Puiseux Algorithm for triple points }\label{NP for triple} 
In this section we use the Newton-Puiseux algorithm to study triple points. 
Let $\bC:f=0$  be a plane curve of degree $d$ and $P$ a triple point for $\bC$.  We focus on the most interesting case, i.e. when $\bC$ has a unique triple tangent $r$ (for an  analysis modulo analytic equivalence of the branches at a triple point see [HH]). }It is not restrictive to assume that $P=O$ and $r: y=0$ , hence
$$f=y^3+\sum_{i+j=4}^d a_{ij}x^iy^j.$$
The cases where $f$ is divisible by $y$ are trivial, so we assume that $y\nmid f$, that is, $f$ is Newton convenient. \b We recall that, by Remark \ref{mult}, the curve $\bC$ can have three 1-branches, or one 2-branch and one 1-branch, or a 3-branch at $O$. If $\bC$ has a 3-branch at $O$, then its parameterization will be of type $(T^3, \sum_{n\geq 4} a_n T^n)$; this follows by Definition \ref{tangent}, since the unique tangent to the branch is $y=0$.

\dd We use the notations of subsection \ref{paths}, with the unique exception that the $g_n$ associated to $f_n^{\textbf{a}_n}$ will be denoted by $g_n^{\textbf{a}_n}$, in order to specify the edge. If $\textbf{a}_n$ is an edge of the Newton polygon, the height of $\textbf{a}_n$ is denoted by $h(\textbf{a}_n)$; the height of the Newton polygon, denoted by $h(\Gamma(f_n))$, is defined to be the sum of the heights of its edges.

\dd We indicate with $_na_{ij}$ the $x^iy_n^j$ coefficient in $f_n$. We also denote by $A_{c_n,r_n}$ the sub-algorithm of Algorithm \ref{Algorithm1} starting inside step $n$ with the choice of the root $q_n= c_nx^{r_n}$ for $g_n$.

\dd In order to describe the output of the algorithm we have to answer the following questions: 
\begin{enumerate}
	\item Which are the possible choices of $c_n$ and $r_n$ in each case?
	\item When, for each choice of $(c_n,r_n)$, can we stop the algorithm $A_{c_n,r_n}$?
\end{enumerate}
The Stop criterion \ref{stop2} is an answer to question 2, but we can improve on it:
\begin{remark}\label{stop}\rm  By Lemma \ref{lemma2}, if $c_n$ has multiciplity 1 as root of $g_n$ then $f_{n+1}$ satisfy the hypothesis of Theorem \ref{ift}, because in $f_{n+1}$ there is the monomial $c_ny_{n+1}$. Hence, if $c_n$ has multiplicity 1 as root of $g_n$, we can stop $A_{c_n,r_n}$ at the step $n-1$, i.e. without calculating $f_{n+1}$. From now on we use this new stop criterion because it is more efficient.
\end{remark}
To answer question 1 we observe that the following statements hold:
\begin{enumerate}[leftmargin=*,label=(\roman*)]
	\item Since the height of Newton Polygon decreases at each step, i.e. $h(\Gamma(f_{n+1})) \leq h(\Gamma(f_n))$, and $h(\Gamma(f))=3$, we have to consider Newton Polygons of height $\leq 3$ only. 
	\item In general $f_n\in\C\left\lbrace \left\lbrace x\right\rbrace \right\rbrace [y]$. However, since $f\in\C[x,y]$, there exists $\overline{n}\in\N$ such that \linebreak $f_n\in\C[x,y]\;\;\forall\;\; n\leq\overline{n}$. Hence for $n\leq \bar n$ the vertices of $\Gamma(f_n)$ are in $\N^2$. 
\end{enumerate}
We start answering question 1 just for $n\leq\bar n$; we shall see later that this is enough. 
\begin{remark}\label{NP cases} \rm The possible shapes of the Newton Polygons arising in the algorithm for a triple point with a triple tangent  and for $n\leq\overline{n}$ are described in the following 11 figures, and each case is discussed below. Namely, in each case we write explicitly the $g_n$ associated to each edge, its roots $c_nx^{r_n}$ and relative multiplicities, so that we have a complete answer to question 1, and then we examine the algorithm $A_{c_n,r_n}$.
\begin{center}
\scalebox{.8}{\includegraphics{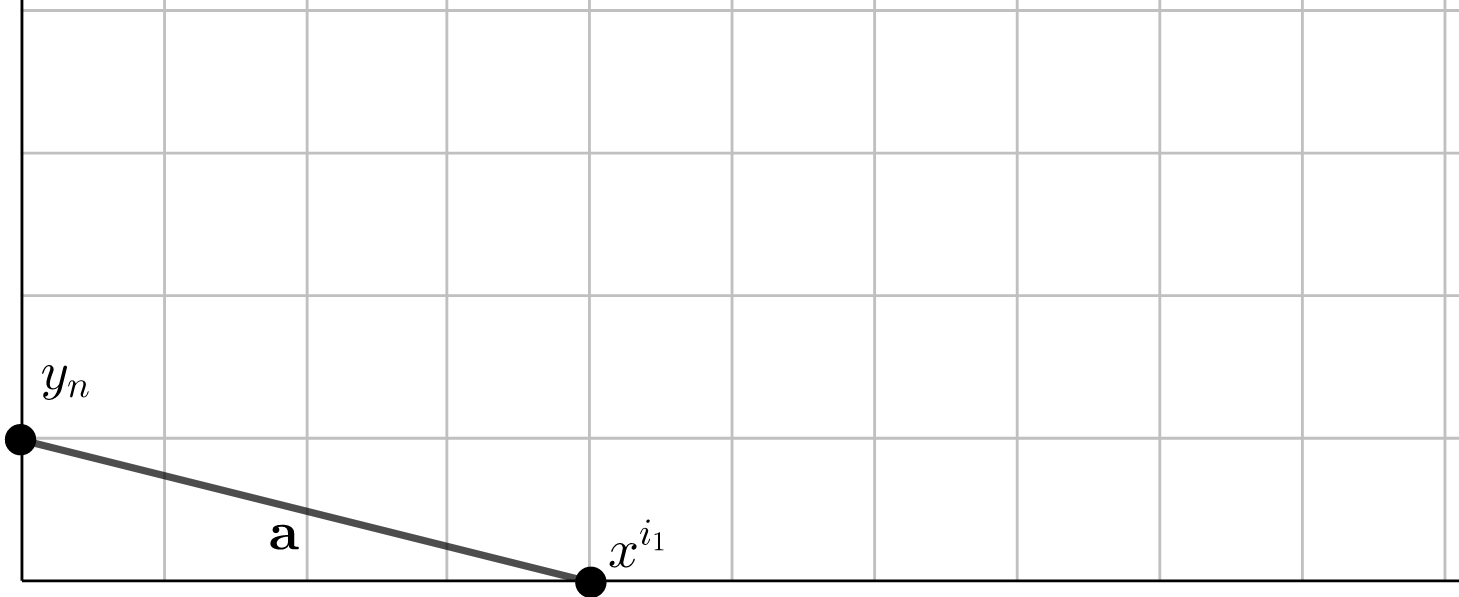}}
\scalebox{.8}{\includegraphics{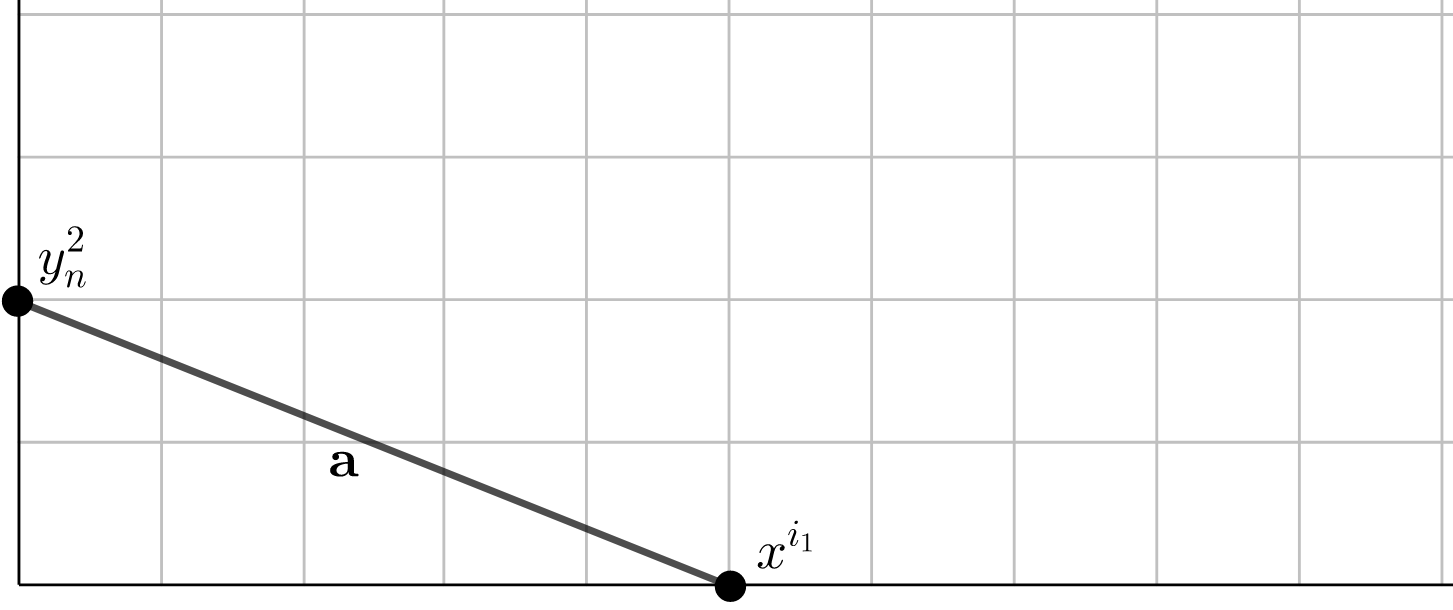}}
\scalebox{.8}{\includegraphics{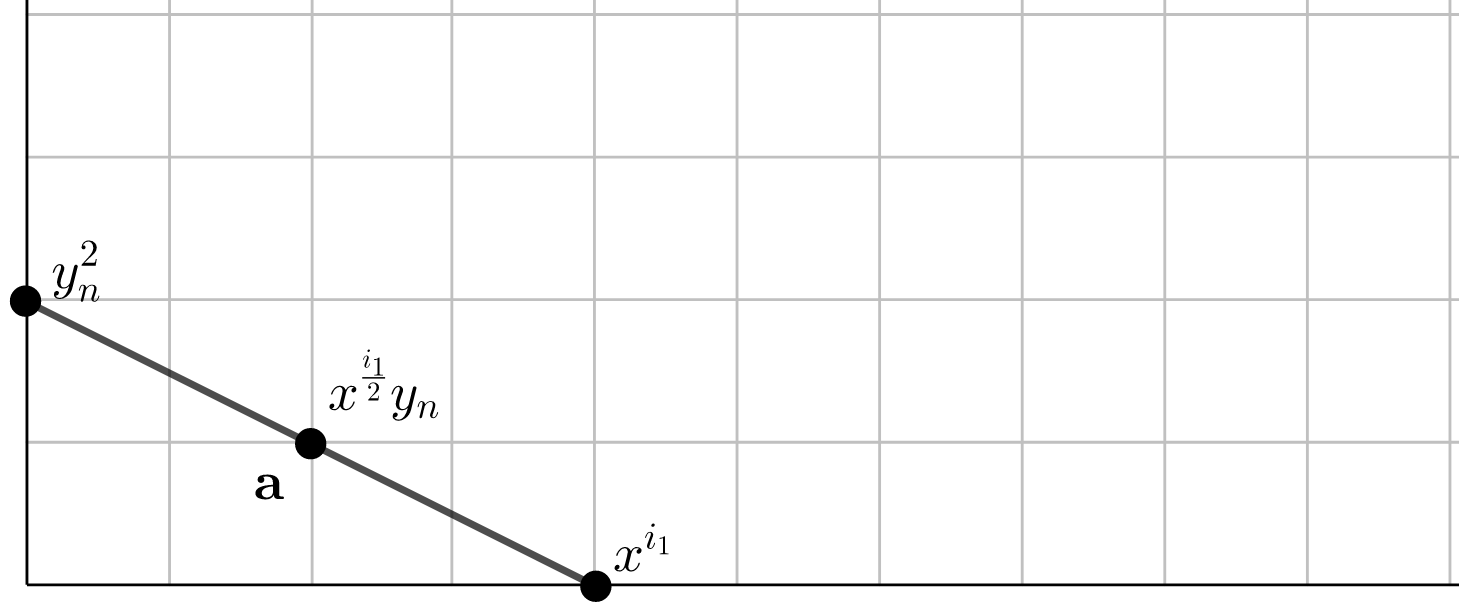}}
\scalebox{.8}{\includegraphics{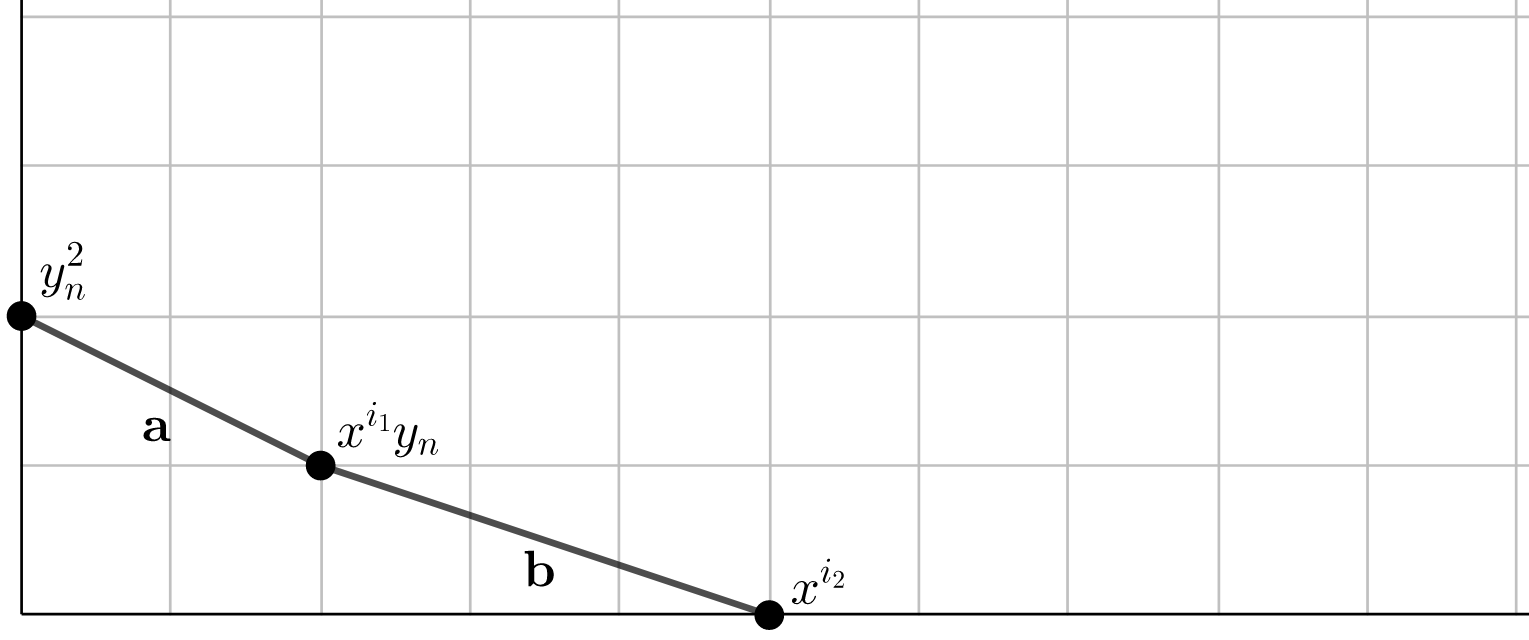}}
\scalebox{.8}{\includegraphics{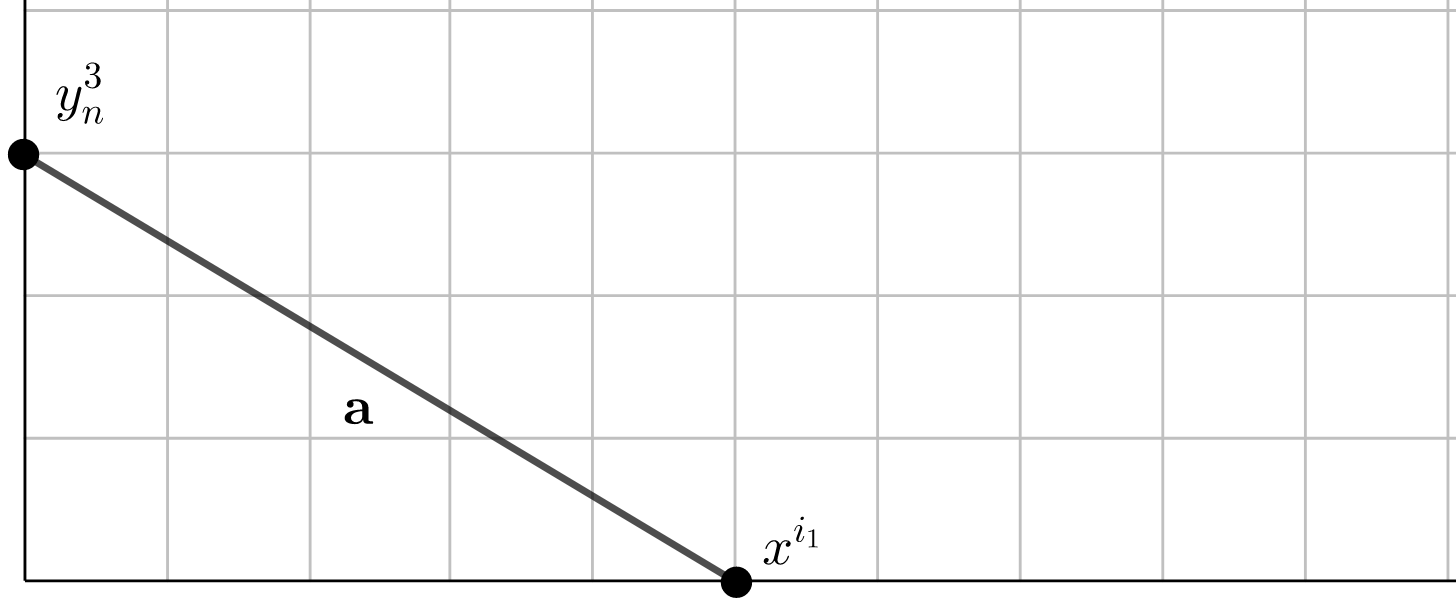}}
\scalebox{.8}{\includegraphics{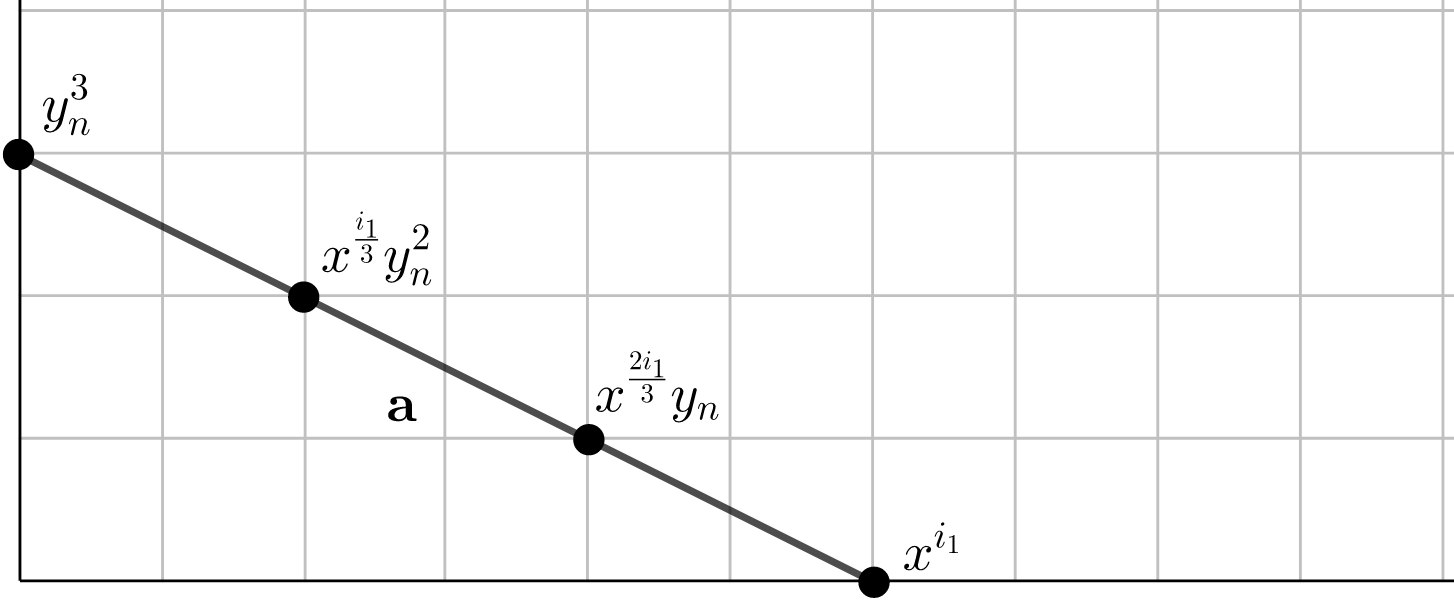}}
\scalebox{.8}{\includegraphics{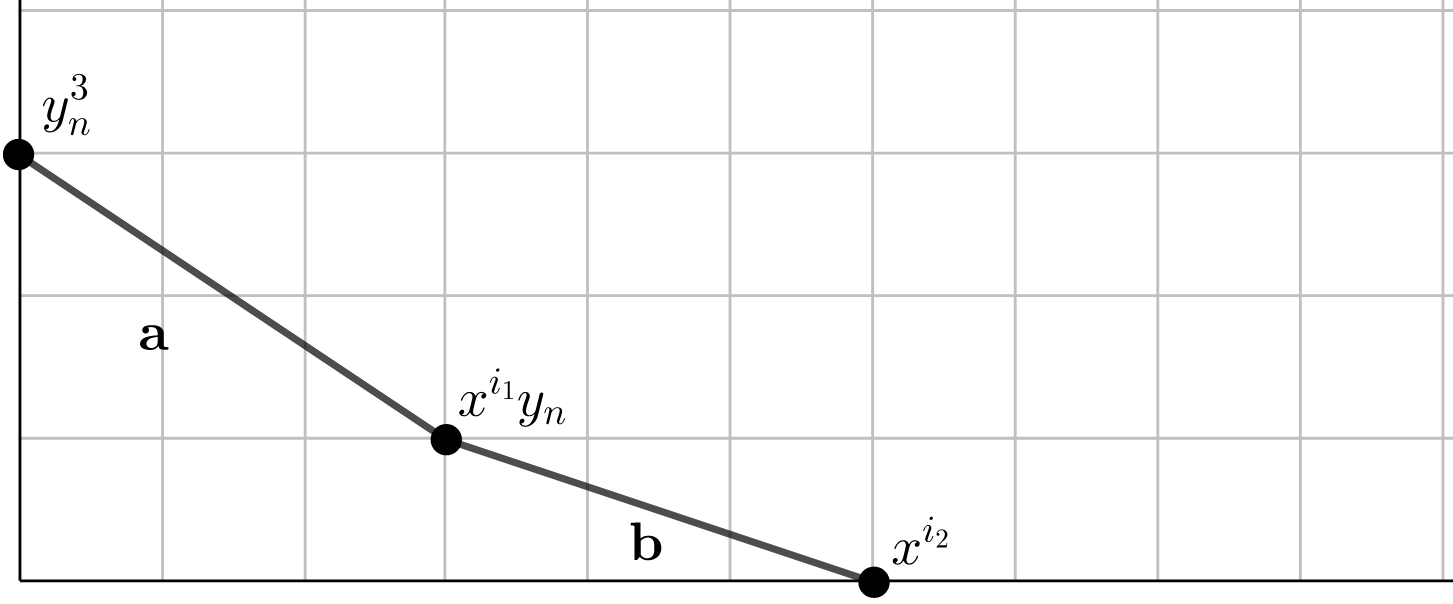}}
\scalebox{.8}{\includegraphics{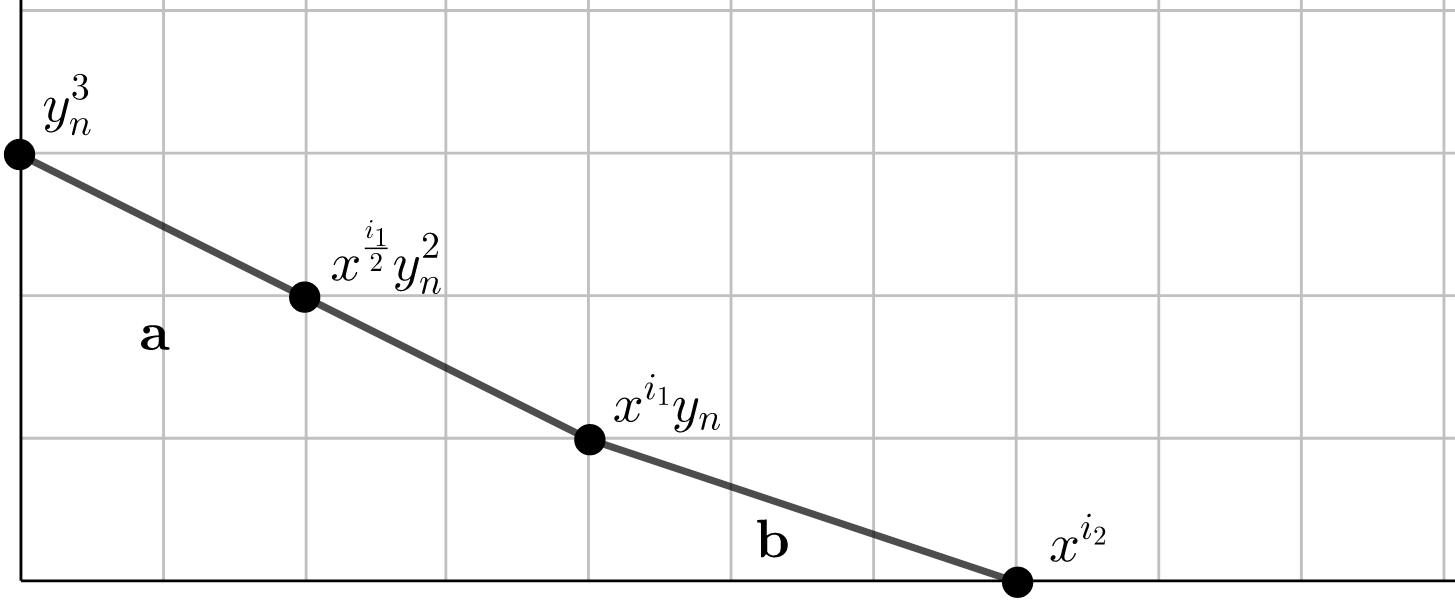}}
\scalebox{.8}{\includegraphics{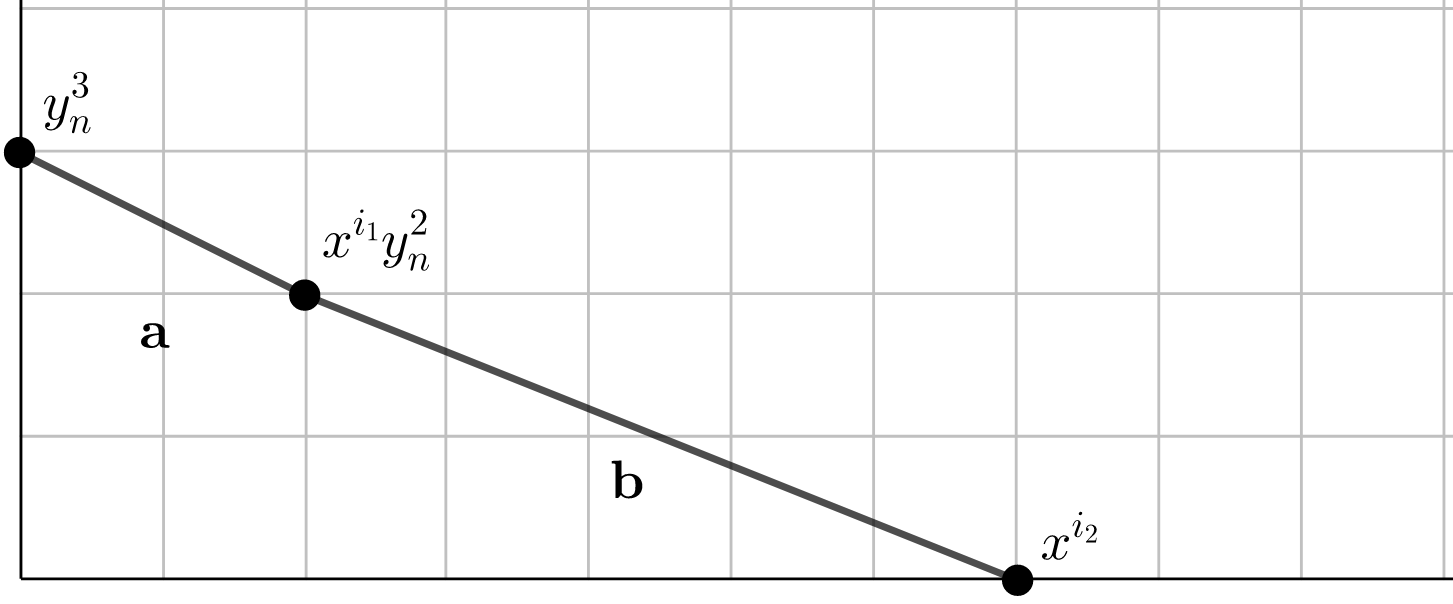}}
\scalebox{.8}{\includegraphics{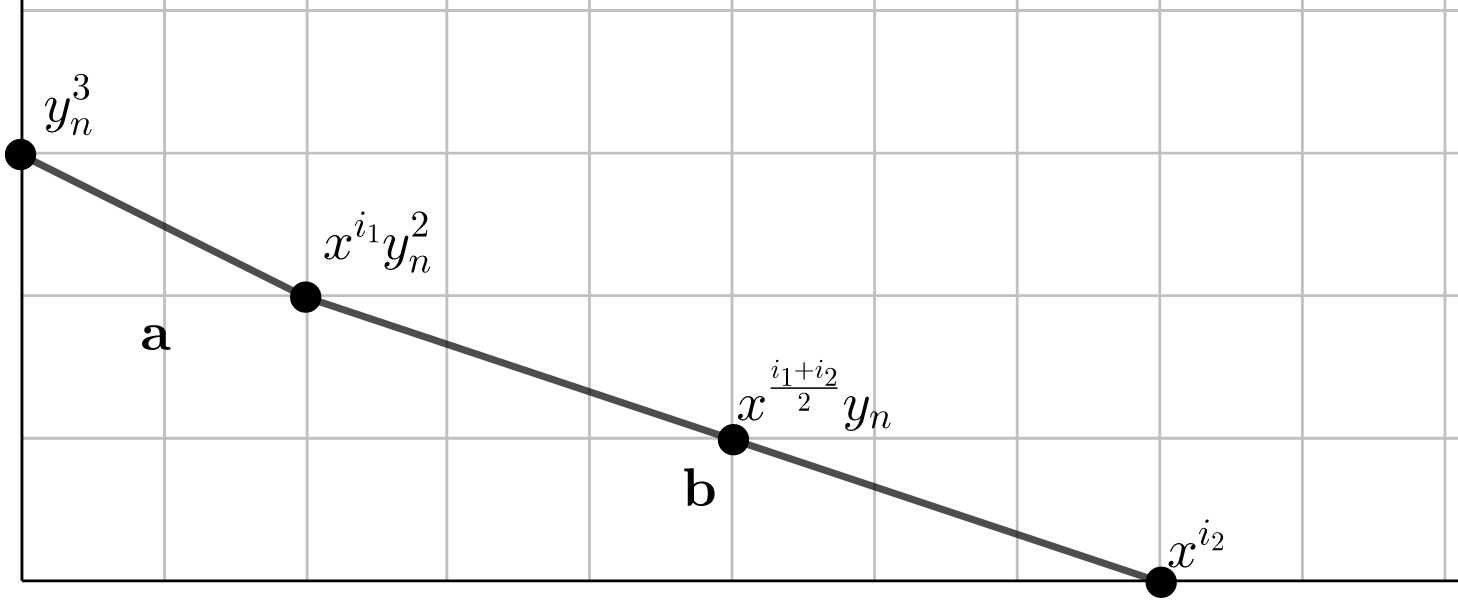}}
\scalebox{.8}{\includegraphics{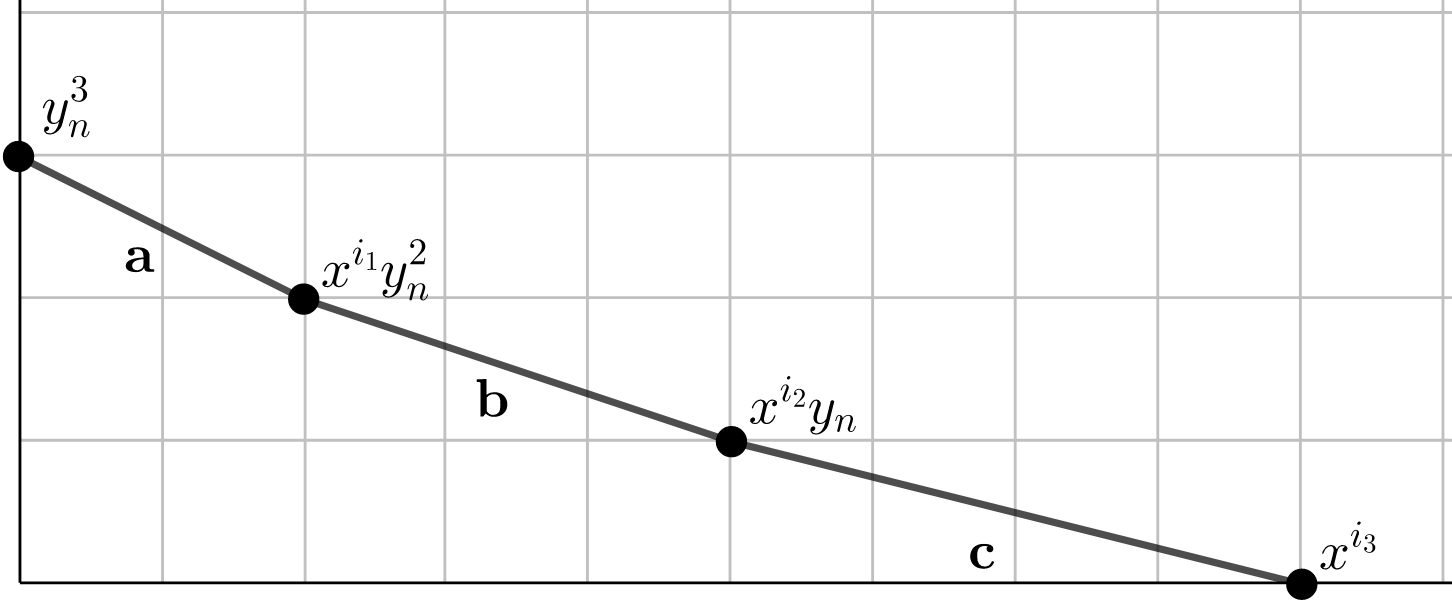}}
\end{center}
\eatit{
\begin{figure}
\centering
\includegraphics[width=.49\textwidth]{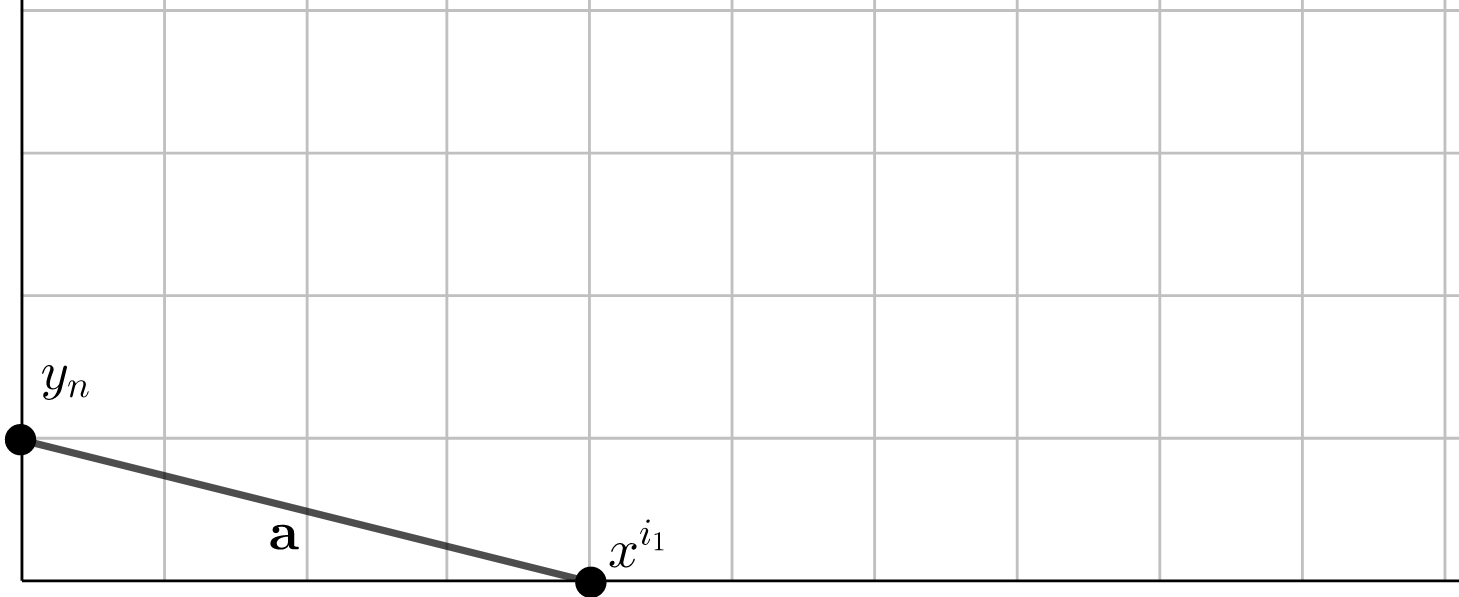}\hfil
\includegraphics[width=.49\textwidth]{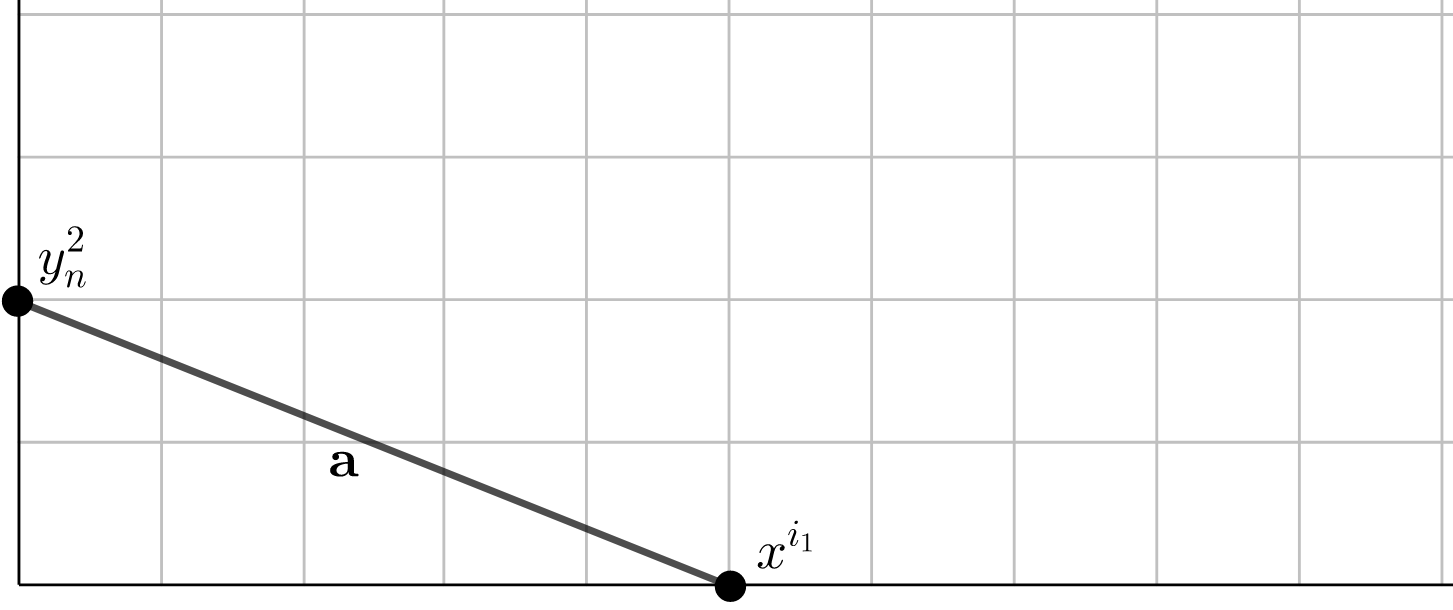}
\caption{Case 1. $(u,v)=(1,1)$ \quad \quad Case 2.1. $(u,v)=(2,1)$ and $i_1\not\equiv 0\Mod{2}$}\label{etichetta}
\end{figure}

\begin{figure}
\centering
\includegraphics[width=.49\textwidth]{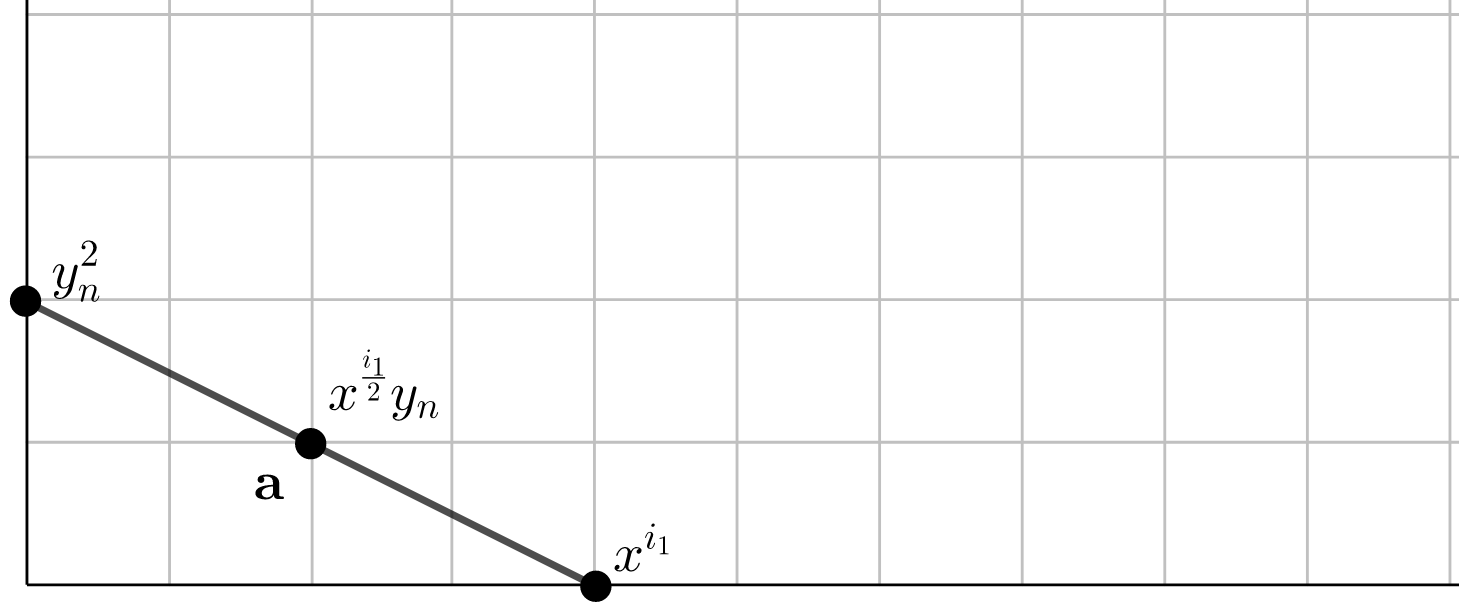}\hfil
\includegraphics[width=.49\textwidth]{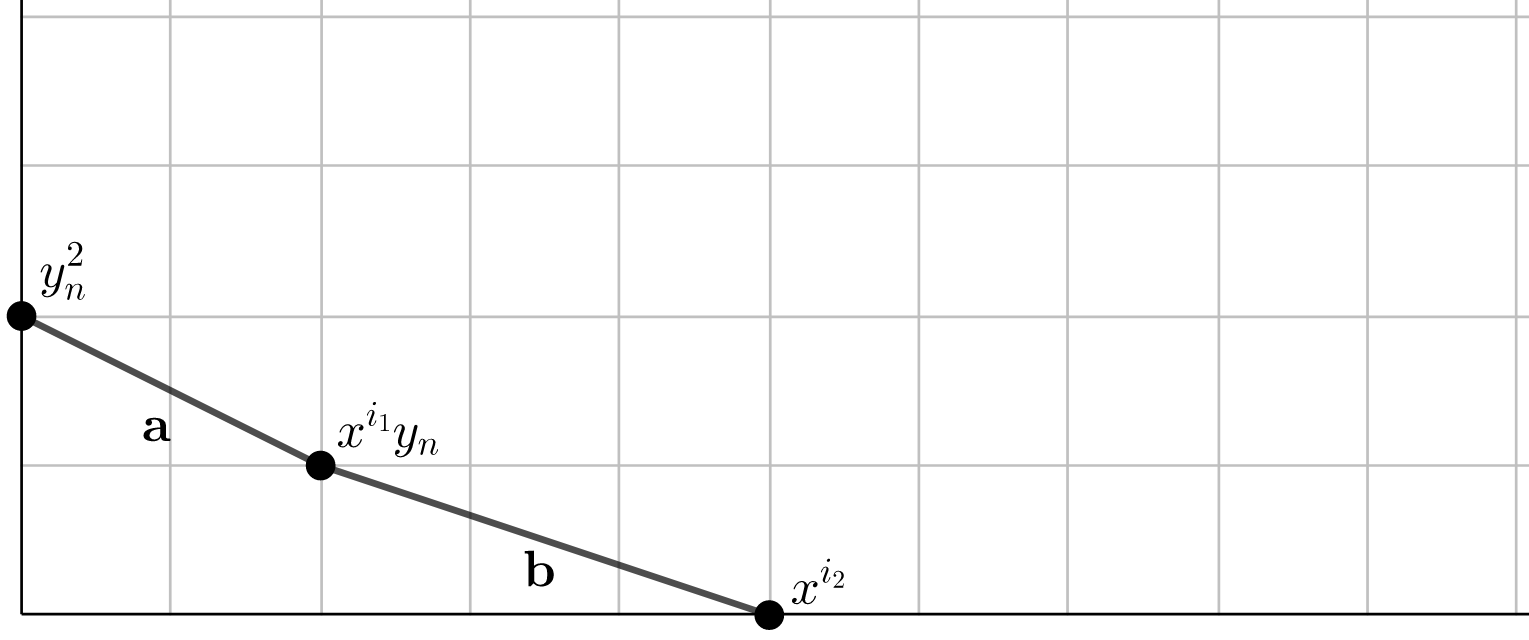}
\caption{Case 2.2. $(u,v)=(2,1)$ and $i_1\equiv 0\Mod{2}$ \quad \quad Case 3. $(u,v)=(2,2)$}\label{}
\end{figure}

\begin{figure}
\centering
\includegraphics[width=.49\textwidth]{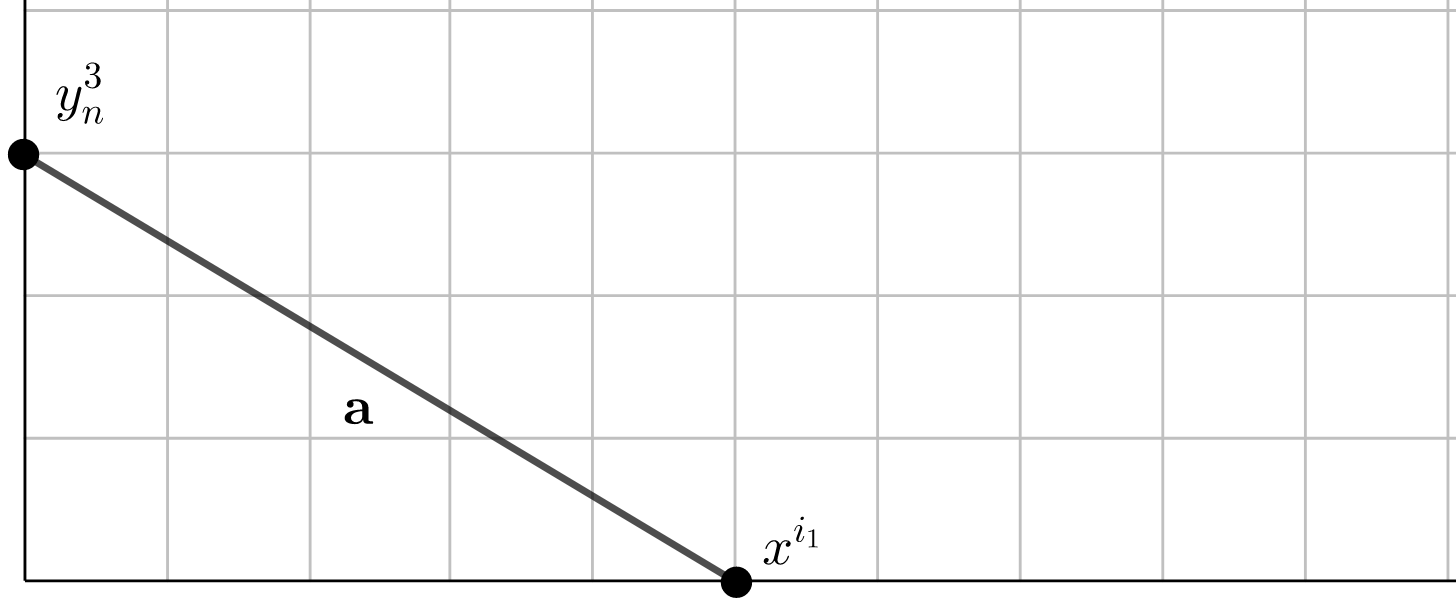}\hfil
\includegraphics[width=.49\textwidth]{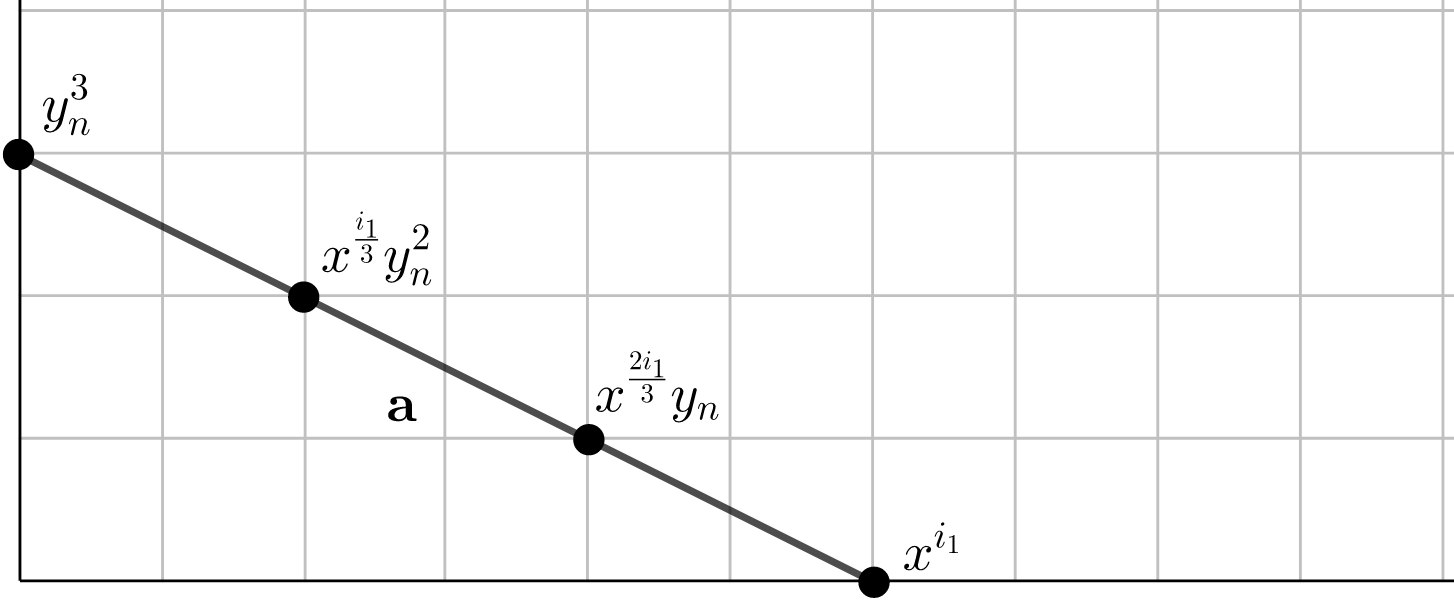}
\caption{Case 4.1. $(u,v)=(3,1)$ and $i_1\not\equiv0\Mod{3}$ \quad \quad Case 4.2.$(u,v)=(3.1)$ and $i_1\equiv0\Mod{3}$}\label{}
\end{figure}

\begin{figure}
\centering
\includegraphics[width=.49\textwidth]{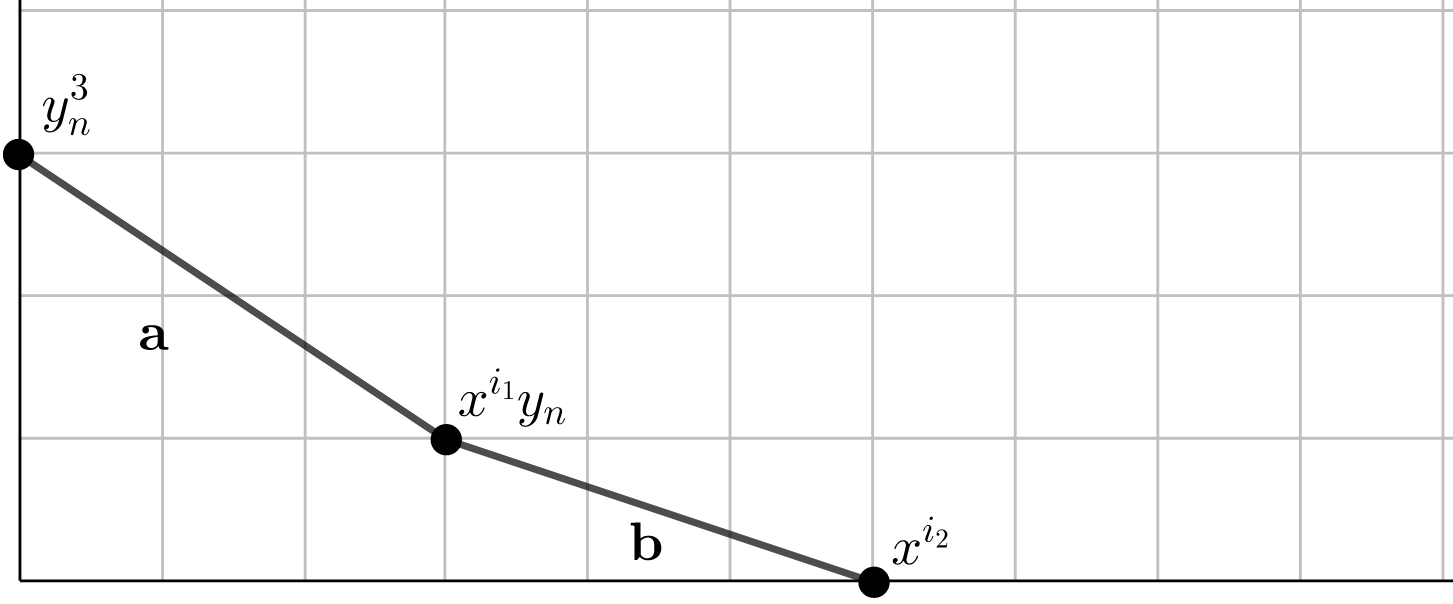}\hfil
\includegraphics[width=.49\textwidth]{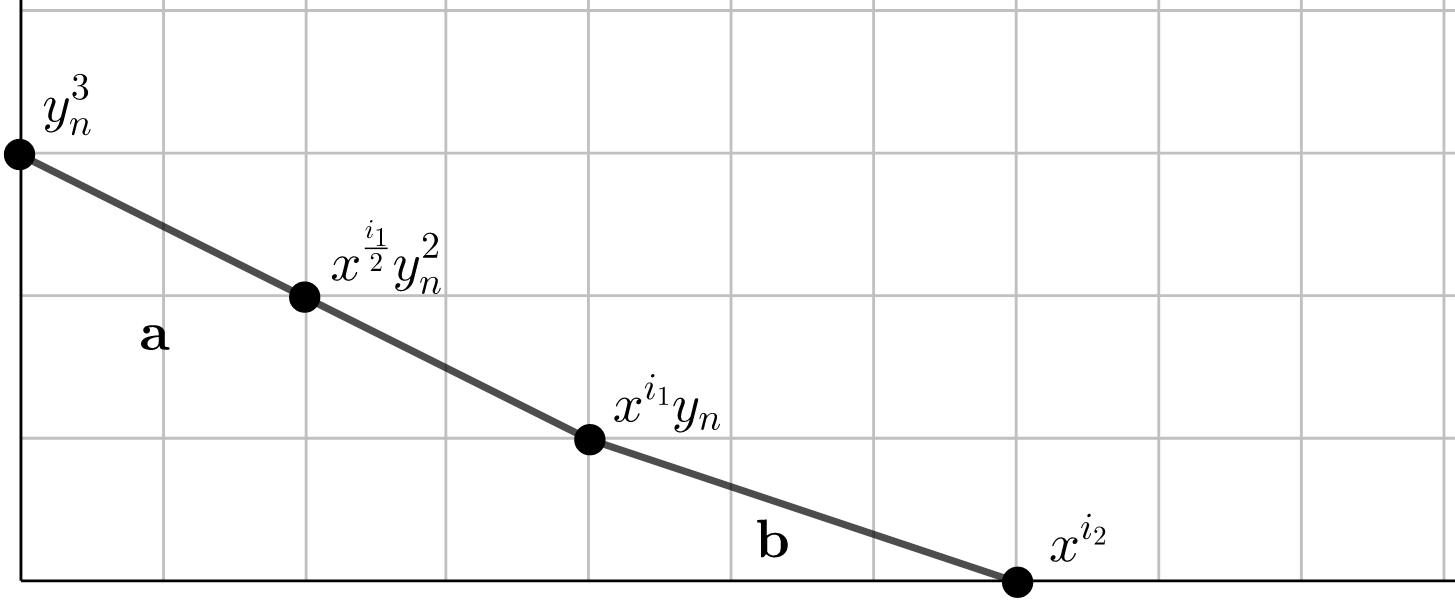}
\caption{Case 5.1. $u=3,v=2(u,v)=(3,2)$,  $h(\textbf{a})= 2, \; i_1\not\equiv0\Mod{2}$ \quad \quad Case 5.2. $(u,v)=(3,2)$, $h(\textbf{a})= 2, \;i_1\equiv0\Mod{2}$}
\end{figure}

\begin{figure}
\centering
\includegraphics[width=.49\textwidth]{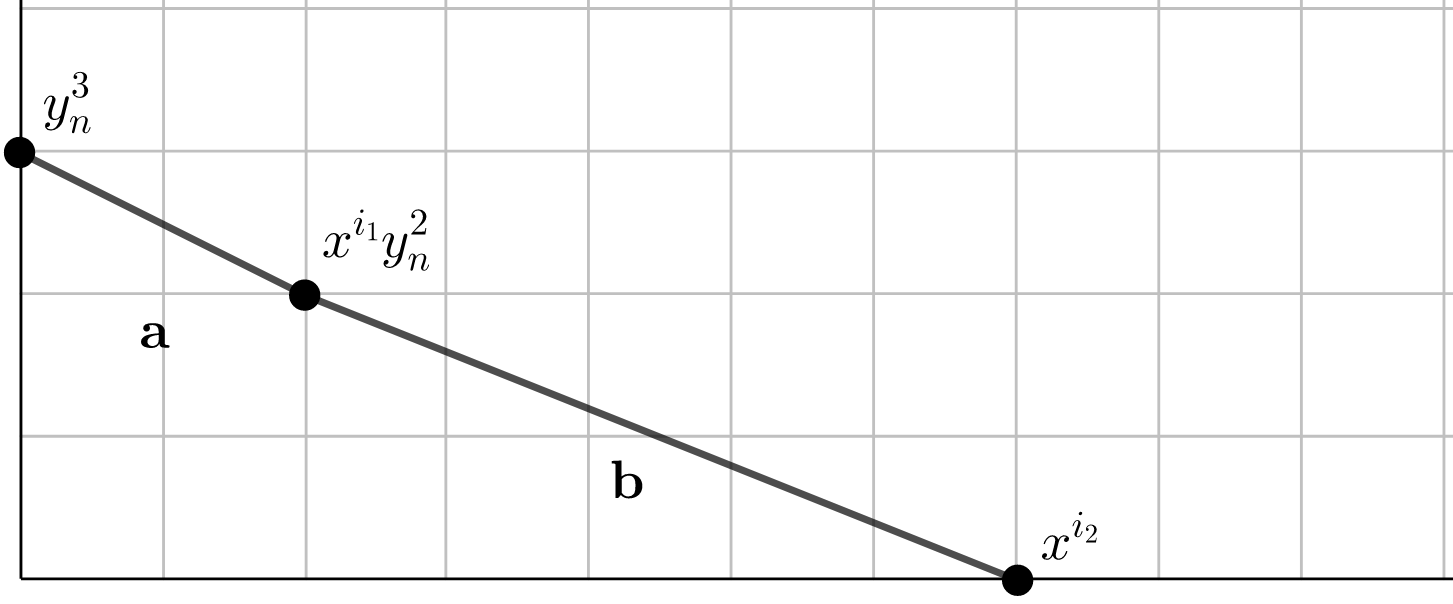}\hfil
\includegraphics[width=.49\textwidth]{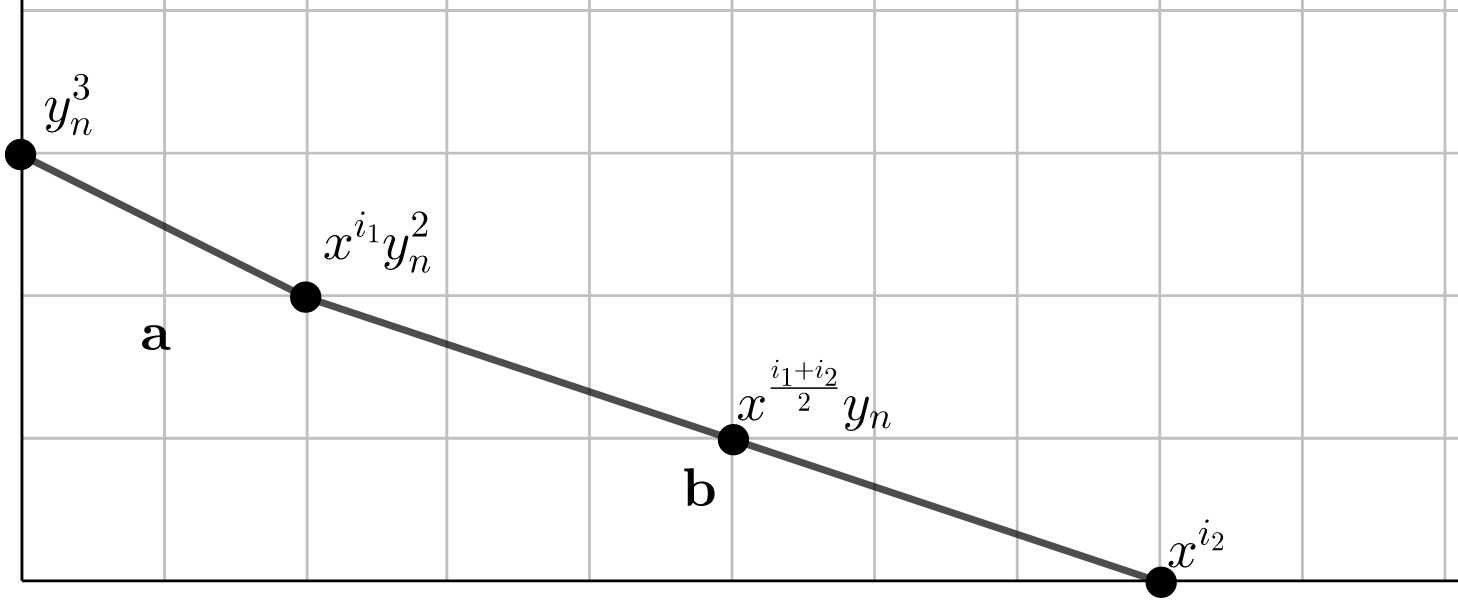}
\caption{Case 6.1. $(u,v)=(3,2)$, the rightmost edge is 2 high and $i_1+i_2\not\equiv0\Mod{2}$ \quad \quad Case 6.2. $(u,v)=(3,2)$, the rightmost edge is 2 high and $i_1+i_2\equiv0\Mod{2}$}\label{}
\end{figure}

\begin{figure}[h!]
	\centering
	\includegraphics[width=.49\textwidth]{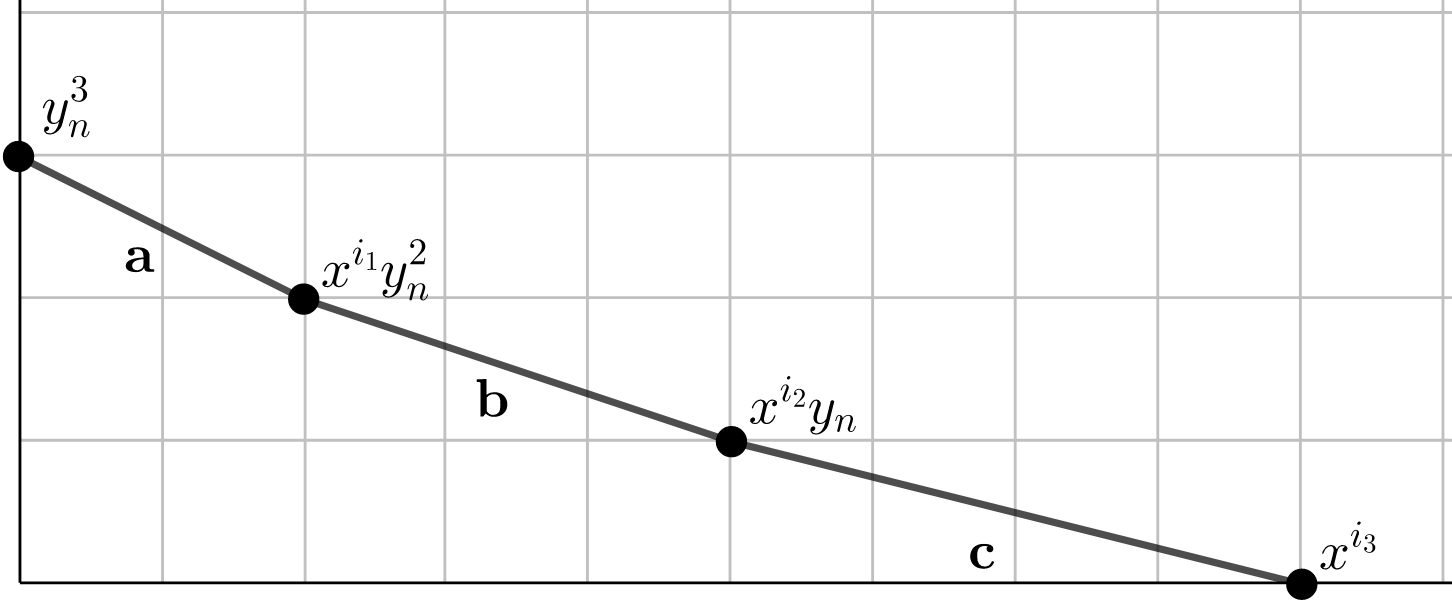}
	\caption{Case 7. $(u,v)=(3,3)$}
\end{figure}

}
%\pagebreak
\bigskip
\textbf{Case 1}: $\Gamma(f_n)= \{\textbf{a}\}$, $h(\textbf{a})=1$, and we have
$$g_{n}^{\textbf{a}}={_na_{01}}y_n+ {_na}_{i_10}x^{i_1}$$
$g_n^{\textbf{a}}$ has a unique root of multiplicity 1,  $c_1x^{i_1}$ with $c_1=-\frac{{_na}_{i_10}}{_na_{01}}$, which by Remark \ref{stop} gives a 1-branch of $\bC$, and the algorithm $A_{c_1,i_1}$ can stop here.

\dd \a \textbf{Case 2.1}:  $\Gamma(f_n)= \{\textbf{a}\}$, $h(\textbf{a})=2$, and $i_1\not\equiv 0\Mod{2}$, so we have

$$g_{n}^{\textbf{a}}={_na_{02}}y_n^2+_n\!a_{i_10}x^{i_1}$$
so that $g_n^{\textbf{a}}$ has two distinct roots of multiplicity 1
$$ c_1x^{\frac{i_1}{2}},\quad c_2x^{\frac{i_1}{2}}\qquad\text{ where } c_1^2=-\frac{{_n}a_{i_10}}{_na_{02}}, \;\, c_2=-c_1.$$
 By Remark \ref{stop}, the algorithms $A_{c_1,\frac{i_1}{2}}$ and $A_{-c_1,\frac{i_1}{2}}$ can stop here; moreover they give equivalent parameterizations of the same 2-branch of $\bC$, since if not, the sum of the multiplicities of the branches at $O$ would be $\geq 4$. 
 The parameterization $(T^2, \dots + c_1T^{i_1}+ \dots)$ composed with the biolomorphic function $T\mapsto -T$ is equal to $(T^2, \dots + c_2T^{i_1}+ \dots)$ ; in fact, since they are equivalent, by Remark \ref{equiv} there is an $\omega$, $\omega^2=1$, such that $R_2$ is obtained from $R_1$ by substituting $T\rightarrow \omega T$; since $c_1(- T)^{i_1}=c_2 T^{i_1}$ we get $\omega=-1$.

\dd \a \textbf{Case 2.2.} $\Gamma(f_n)= \{\textbf{a}\}$,  $h(\textbf{a})=2$, and $i_1\equiv 0\Mod{2}$, so we have

$$g_{n}^{\textbf{a}}={_na_{02}}y_n^2+{_na_{\frac{i_1}{2}1}}x^{\frac{i_1}{2}}y_n+{_na_{i_10}}x^{i_1}.$$
Let $c_1x^{i_1\over 2}, c_2x^{i_1\over 2}$ denote the roots of $g_{n}^{\textbf{a}}$; we distinguish two cases:

	\bigskip \textbf{Case 2.2.1.} If $c_1\neq c_2$, i.e. if $g_{n}^{\textbf{a}}$ has two roots of multiplicity 1, they give 2 distinct 1-branches of $\bC$
	and the algorithms $A_{c_1,\frac{i_1}{2}}$ and $A_{c_2,\frac{i_1}{2}}$ can stop here.
	\b Notice that if ${_na_{\frac{i_1}{2}1}}=0$ we apparently are in the same situation as {Case 2.1}, but here $k:={i_1\over 2}$ is a positive integer, the two simple roots are
$$c_1x^k,\quad -c_1x^k$$
and the biolomorphic function $T\mapsto -T$ used in 2.1 doesn't work here: it does not give the equivalence of $(T, \dots+c_1T^{k}+ \dots)$ and $(T, \dots-c_1T^{k}+ \dots)$. In fact, the two simple roots correspond to two not equivalent parameterizations of a 1-branch each.

	\bigskip \textbf{Case 2.2.2.} If $c_1=c_2$, i.e. if $g_{n}^{\textbf{a}}$ has one double root, the algorithm $A_{c_1,\frac{i_1}{2}}$ has to go on:
	it can give either two 1-branches or one 2-branch of $\bC$, hence globally we may find either three  1-branches or one 2-branch and one 1-branch.
Notice that $f_{n+1}$ is still in $\C[x,y]$, since $i_1\equiv 0\Mod{2}$.

\bigskip
\noindent \textbf{Case 3.} $\Gamma(f_n) = \{\textbf{a}, \textbf{b}\}$, $h(\textbf{a})=1, h(\textbf{b})=1$; we have

$$g_n^{\textbf{a}}={_na_{02}}y_n+{_na}_{i_11}x^{i_1}, \qquad  g_n^{\textbf{b}}={_na}_{i_11}y_n+{_na_{i_20}}x^{i_2-i_1}.$$
Both $g_n^{\textbf{a}}$ and $g_n^{\textbf{b}}$ have a unique root of multiplicity 1,  $-\frac{_na_{i_11}}{_na_{02}}x^{i_1}$ and $-\frac{_na_{i_20}}{_na_{i_11}}x^{i_2-i_1}$, which give a 1-branch of $\bC$ each, and the algorithms can stop here.

\dd \a \textbf{Case 4.1.} $\Gamma(f_n)=\{\textbf{a}\}$, $h(\textbf{a})=3$, and $i_1\not\equiv 0\Mod{3}$, so we have
$$g_n^{\textbf{a}}={_na_{03}}y_n^3+{_na_{i_10}}x^{i_1}$$
$g_n^{\textbf{a}}$ has three distinct roots of multiplicity 1, namely:
$$c_1x^{\frac{i_1}{3}},\quad c_2x^{\frac{i_1}{3}},\quad c_3x^{\frac{i_1}{3}}\qquad \text{where } c_1^3=-\frac{_na_{i_10}}{_na_{03}}, \;\;c_2=e^{i{2\pi\over 3}} c_1, \;\; c_3=e^{i{4\pi\over 3}} c_1$$
and analogously to Case 2.1, all of them give equivalent parameterizations of the same 3-branch of $\bC$, since if not, the sum of the multiplicities of the branches at $O$ would be $\geq 6$; for example the parameterization $(T^3, \dots + c_1T^{i_1}+ \dots)$ composed with the biolomorphic function $T\mapsto e^{i{2\pi\over 3}}T$ is equal to $(T^3, \dots + c_2T^{i_1}+ \dots)$.
Since the three roots have multiplicity 1, we can stop the algorithms.

\dd \a \textbf{Case 4.2.} $\Gamma(f_n)=\{\textbf{a}\}$, $h(\textbf{a})=3$, and $i_1\equiv 0\Mod{3}$, so we have
$$g_n^\textbf{a}={_na_{03}}y_n^3+{_na}_{\frac{i_1}{3}2}x^{\frac{i_1}{3}}y_n^2+{_na_{\frac{2i_1}{3}1}}x^{2\frac{i_1}{3}}y_n+{_na_{i_10}}x^{i_1}.$$

\b Let $c_1x^{i_1\over 3},\quad c_2x^{i_1\over 3}, \quad c_3x^{i_1\over 3}$ denote the roots of $g_{n}^{\textbf{a}}$; we distinguish three cases:

\bigskip \textbf{Case 4.2.1.}  If $c_1\neq c_2 \neq c_3\neq c_1$,  i.e. if $g_{n}^{\textbf{a}}$ has three roots of multiplicity 1, they give 3 distinct 1-branches of $\bC$ so the algorithms can stop here.
\b Notice that if ${_na}_{\frac{i_1}{3}2}={_na_{\frac{2i_1}{3}1}} =0$, we apparently are in the same situation as {Case 4.1}, but here $k:={i_1\over 3}$ is a positive integer, the three simple roots are
$$c_1x^k,\quad e^{i{2\pi\over 3}} c_1x^k,\quad e^{i{2\pi\over 3}} c_1x^k$$
and the biolomorphic function $T\mapsto e^{i{2\pi\over 3}}T$ used in 4.1 doesn't work here: it does not give the equivalence of $(T,\dots+ c_1T^{k}+ \dots)$ and $(T, \dots+c_2T^{k}+ \dots)$. In fact, the three simple roots correspond to three not equivalent parameterizations of a 1-branch each.

	\bigskip \textbf{Case 4.2.2.}  If $c_1\neq c_2=c_3$, ,  i.e. if $g_{n}^{\textbf{a}}$ has one double and one simple root, we can stop the algorithm $A_{c_1,\frac{i_1}{3}}$, which gives a 1-branch of $\bC$. 
	On the other hand, $A_{c_2,\frac{i_1}{3}}$ cannot be stopped because the root has multiplicity 2; going on, it can give two 1-branches or one 2-branch of $\bC$. Notice that in this case $f_{n+1}$ is still in $\C[x,y]$, since $i_1\equiv 0\Mod{3}$.

	\bigskip \textbf{Case 4.2.3.} 	
 If $c_1=c_2=c_3$,  i.e. if $g_{n}^{\textbf{a}}$ has one triple root, the algorithm $A_{c_1,\frac{i_1}{3}}$ has to go on. Notice that in this case $f_{n+1}$ is still in $\C[x,y]$, since $i_1\equiv 0\Mod{3}$.

\dd \a  \textbf{Case 5.1.} 
$\Gamma(f_n) = \{\textbf{a}, \textbf{b}\}$, $h(\textbf{a})=2, h(\textbf{b})=1$ and $i_1\not\equiv 0\Mod{2}$, so we have
$$g_n^\textbf{a}={_na_{03}}y_n^2+{_na_{i_11}}x^{i_1}, \qquad g_n^\textbf{b}={_na_{i_11}}y_n+{_na_{i_20}}x^{i_2-i_1}$$
$g_n^{\textbf{a}}$ has two distinct roots of multiplicity 1: $c_1x^{i_1\over 2}, c_2x^{i_1\over 2}$ with $c_1^2=-\frac{_na_{i_11}}{_na_{03}},\;\; c_2=-c_1$, which give equivalent parameterizations of the same 2-branch of $\bC$, while $g_n^{\textbf{b}}$ has one root of multiplicity 1: 
$ -\frac{_na_{i_20}}{_na_{i_11}}x^{i_2-i_1}$, giving one 1-branch. The algorithms can stop here.

\dd \a \textbf{Case 5.2.} $\Gamma(f_n) = \{\textbf{a}, \textbf{b}\}$,  $h(\textbf{a})=2, h(\textbf{b})=1$ and $i_1\equiv 0\Mod{2}$, so we have
$$g_n^\textbf{a}={_na_{03}}y_n^2+{_na_{\frac{i_1}{2}2}}x^{\frac{i_1}{2}}y_n+{_na_{i_11}}x^{i_1}, \qquad g_n^\textbf{b}={_na_{i_11}}y_n+{_na_{i_20}}x^{i_2-i_1}$$
$g_n^{\textbf{b}}$ has one root of multiplicity 1, which  gives a 1-branch of $\bC$ and the relative algorithm can stop here, while,   if $c_1x^{i_1\over 2}, c_2x^{i_1\over 2}$ denote the roots of $g_{n}^{\textbf{a}}$, we distinguish two cases (notice that ${_na_{\frac{i_1}{2}2}}$ may be 0, in which case we are apparently in the same situation as Case 5.1):

	\bigskip\textbf{Case 5.2.1.} If $c_1\neq c_2$, i.e. if $g_{n}^{\textbf{a}}$ has two roots of multiplicity 1, each of them gives a 1-branch of $\bC$
	so we have globally three 1-branches of $\bC$ and the algorithms can stop here.
	
	\bigskip \textbf{Case 5.2.2.} If $c_1=c_2$, i.e. if $g_{n}^{\textbf{a}}$ has the double root $c_1x^{\frac{i_1}{2}}$, the algorithm $A_{c_1,\frac{i_1}{2}}$ can't stop; going on, it can give two 1-branches or one 2-branch of $\bC$, hence globally we may find three 1-branches or one 2-branch and one 1-branch.
Notice that in this case $f_{n+1}$ is still in $\C[x,y]$, since $i_1\equiv 0\Mod{2}$.

\dd \a \textbf{Case 6.1.} $\Gamma(f_n) = \{\textbf{a}, \textbf{b}\}$,   $h(\textbf{a})=1, h(\textbf{b})=2$ and $i_1+i_2\not\equiv 0\Mod{2}$, so we have
$$g_n^\textbf{a}={_na_{03}}y_n+{_na}_{i_12}x^{i_1}, \qquad g_n^\textbf{b}={_na_{i_12}}y_n^2+{_na_{i_20}}x^{i_2-i_1}.$$
Since $i_2-i_1\not\equiv 0\Mod{2}$ the discussion is analogous to Case 5.1: $\bC$ has one 1-branch and one 2-branch, and the algorithms can stop here. 

\dd \a \textbf{Case 6.2.} $\Gamma(f_n) = \{\textbf{a}, \textbf{b}\}$, $h(\textbf{a})=1, h(\textbf{b})=2$ and $i_1+i_2\equiv 0\Mod{2}$, so we have
$$g_n^\textbf{a}={_na_{03}}y_n+{_na}_{i_12}x^{i_1}, \qquad g_n^\textbf{b}={_na_{i_12}}y_n^2+{_na_{\frac{i_1+i_2}{2}1}}x^{\frac{i_2-i_1}{2}}y_n+{_na_{i_20}}x^{i_2-i_1}.$$
Notice that the assumption $i_1+i_2\equiv 0\Mod{2}$ gives $i_2-i_1\equiv 0\Mod{2}$, hence the discussion is analogous to Case 5.2, namely we distinguish two cases:

\bigskip  \textbf{Case 6.2.1.} if $g_{n}^{\textbf{b}}$ has two roots of multiplicity 1, each of them gives a 1-branch of $\bC$
so we have globally three 1-branches of $\bC$ and the algorithms can stop here. Notice that ${_na_{\frac{i_1+i_2}{2}1}}$ may be 0, in which case we are apparently in the same situation as Case 6.1, but here $i_2-i_1\equiv 0\Mod{2}$.

\bigskip \textbf{Case 6.2.2.}  if $g_{n}^{\textbf{b}}$ has one double root $c_1x^{\frac{i_1}{2}}$, the relative algorithm can't stop; going on, it can give two 1-branches or one 2-branch of $\bC$, hence globally we may find three 1-branches or one 2-branch and one 1-branch. Notice that in this case $f_{n+1}$ is still in $\C[x,y]$, since $i_1+i_2\equiv 0\Mod{2}$.

\dd \a \textbf{Case 7.} $\Gamma(f_n)=\{\textbf{a}, \textbf{b}, \textbf{c}\}$, $h(\textbf{a})= h(\textbf{b})=h(\textbf{c})=1$ and we have
$$g_n^\textbf{a}={_na_{03}} y_n+{_na_{i_12}}x^{i_1}, \qquad g_n^\textbf{b}={_na_{i_12}}y_n+{_na_{i_21}}x^{i_2-i_1}, \qquad g_n^\textbf{c}={_na_{i_21}}y_n+{_na_{i_30}}x^{i_3-i_2}$$
$g_n^\textbf{a}$, $g_n^\textbf{b}$ and $g_n^\textbf{c}$ have one simple root each and the three roots give a 1-branch of $\bC$ each: the algorithms can stop here.

\end{remark}
Now we are ready to see what happens when we run the Newton-Puiseux Algorithm \ref{Algorithm1} in the case of a triple point with a triple tangent. 
%fine eatit

\begin{thm}\label{tripli} Let $\;\; \bC:   f(x,y)=0\;\;$ be a plane curve of degree $d$, such that $O$ is a triple point for $\bC$ and the $x$-axis is the triple tangent at $O$, so that 
$$f(x,y)=   y^3+\sum_{i+j=4}^d a_{ij}x^iy^j.$$
With the previous notations, the output of the Newton-Puiseux Algorithm \ref{Algorithm1} is as follows:

\dd a finite sequence of $n$ steps like case {\rm \bf 4.2.3} with $n\geq 0$, and then

\dd $(i)$ step $n+1$ is a case {\rm \bf 4.1}; we can stop the algorithm:
\b$\bC$ has one 3-branch at $O$, and a parameterization for it is of the form
$$\left(T^3, b_0T^{3h_0}+\dots+b_nT^{3h_n}+ c_1T^{s}+ \dots \right)\quad s\not\equiv 0\Mod{3}$$

\dd $(ii)$  step $n+1$ is one of the following cases: {\rm \bf 4.2.1}, {\rm \bf 5.2.1}, {\rm \bf 6.2.1}, {\rm \bf 7}; we can stop the algorithm:
\b $\bC$ has three 1-branches at $O$, each of which with a parameterization of the form 
$$\left(T, b_0T^{h_0}+\dots+b_nT^{h_n}+ c_jT^{qj}\dots \right)$$
where $q_1=q_2=q_3$ in case {\rm \bf 4.2.1}, $q_1=q_2$ in cases {\rm \bf 5.2.1}, {\rm \bf 6.2.1}.

\dd $(iii)$  step $n+1$ is one of the following cases: {\rm \bf 5.1}, {\rm \bf 6.1}; we can stop the algorithm:
\b $\bC$ has one 1-branch and one 2-branch at $O$, with parameterizations of the form:
$$\left(T, b_0T^{h_0}+\dots+b_nT^{h_n}+ c_1T^{q}+ \dots \right)$$
 $$\left(T^2, b_0T^{2h_0}+\dots+b_nT^{2h_n}+ d_1T^{2h+1}+\dots \right).$$

\dd $(iv)$ step $n+1$ is one of the following cases: {\rm \bf 4.2.2}, {\rm \bf 5.2.2}, {\rm \bf 6.2.2};  in order to decide if $\bC$ has one 1-branch and one 2-branch or three 1-branches at $O$, the sub-algorithm relative to the double root must go on, so that we find 
\b a finite sequence of $m \geq 0$ steps like case {\rm \bf 2.2.2}, and then:
 \b if the next step is a case {\rm \bf 2.1}, the algorithm stops: $\bC$ has one 1-branch and one 2-branch at $O$, with parameterizations of the form
 $$\left(T, b_0T^{h_0}+\dots+b_nT^{h_n}+ c_1T^{p}+ \dots \right)$$
 $$\left(T^2, b_0T^{2h_0}+\dots+b_nT^{2h_n}+ c_2T^{2p}+d_1T^{2k_1}+\dots+d_mT^{2k_m}+a_1T^{2h+1}+ \dots \right)$$

 \b if the next step is a case {\rm \bf 2.2.1} or {\rm \bf3}, the algorithm stops: $\bC$ has three 1-branches at $O$, with parameterizations of the form
 $$\left(T, b_0T^{h_0}+\dots+b_nT^{h_n}+ c_1T^{p}+ \dots \right)$$
 $$\left(T, b_0T^{h_0}+\dots+b_nT^{h_n}+ c_2T^p+d_1T^{k_1}+\dots+d_mT^{k_m}+e_1T^{q_1}+ \dots \right)$$
 $$\left(T, b_0T^{h_0}+\dots+b_nT^{h_n}+ c_2T^p+ d_1T^{k_1}+\dots+d_mT^{k_m}+g_1T^{q_2}+ \dots \right)$$
where $q_1=q_2$ in case {\rm \bf 2.2.1}.

\end{thm}
 \proof We  claim that, in order to compute the number and multiplicities of $\bC$ at a triple point with a triple tangent, it is enough to consider Newton polygons whose vertices are in $\N^2$; that is, applying the Newton - Puiseux method to $f$, we encounter only the cases described in Remark  \ref{NP cases} above. This is true for $n=0$, since $f\in \C[x,y]$. Now choose a path as described in Algorithm \ref{Algorithm1}, and assume it is true for the steps $1,\dots, n$; hence the Newton Polygon of $f_n$ is one of those described above. If the algorithm can't stop we are in one of the following 5 cases: {\bf 2.2.2, 4.2.2, 4.2.3, 5.2.2, 6.2.2}, and since in each of these cases the vertices of  the Newton Polygon of $f_{n+1}$ are still in $\N^2$, it is true for $n+1$. On the other hand, if the algorithm can stop we don't need to go any further, so the claim is true.
 
 \a We analyze what happens in the assumption that $y_n \nmid f_n$; if, on the contrary, $y_n\mid f_n $, there is a branch coming from an algebraic component of $\bC$, and thus the corresponding series is in fact a polynomial. 

\a Assume that after $n$ steps we are in case {\bf 4.2.3}. Since $3\geq h(\G(f))\geq h(\G(f_m))\geq h(\G(f_n))=3$ for $m\leq n$, we have $h(\G(f_m))=3$ for $m\leq n$, hence each previous step was one of the cases between {\bf 4.1} and {\bf 7}. If one of the previous steps was a case among {\bf 4.1, 4.2.1, 5.1, 5.2.1, 6.1, 6.2.1, 7 } we would have already stopped. If one of the previous steps was a case among {\bf 4.2.2, 5.2.2,  6.2.2},  only the sub-algorithm relative to the double root would have gone on, and in the next steps the height of the Newton Polygon would have been 2. So all the previous steps were of type {\bf 4.2.3}.
\b Let us start with the step 0 which is a case {\bf 4.2.3}, and set $k={i_1\over 3}$; since $i_1\equiv{0}\Mod{3}$, $k$ is a positive integer, and for a suitable coefficient $b$ we have 
$$f= (y-bx^{k})^3+ \sum_{i+kj>3k} a_{ij} x^iy^j.$$
Applying the algorithm we get: $g^\textbf{a}=f^\textbf{a} = (y-bx^k)^3$, $q_0= bx^k$, 
$$  f(x, x^k(b+y_1))= (x^k(b+y_1)-bx^k)^3+ \sum_{i+kj>3k} a_{ij} x^i(x^k(b+y_1))^j=x^{3k}y_1^3+\sum_{i+kj>3k} a_{ij} x^{i+jk}(b+y_1)^j$$
and since each power of $x$ has an exponent $\geq 3k$ we find
$$f_1= {x^{3k}y_1^3+\sum_{i+kj>3k} a_{ij} x^{i+jk}(b+y_1)^j\over x^{3k}}= y_1^3+\dots$$
If  $f_1$ is a case {\bf 4.2.3} again we have:
$$f_1= (y-b_1x^h)^3+ \sum_{i+hj>3h} b_{ij} x^iy^j$$
and we find $q_1= b_1x^{k_1}$ with $k_1 \in \N$; hence if the sequence of cases {\bf 4.2.3} goes on up to step $n$, we have a sequence of roots, setting $b_0:=b, k_0:=k$ (notice that $k_0\geq 2$, otherwise the tangent cone of $C$ at $O$ would be different from $y^3=0$)
 
$$q_0= b_0x^{k_0},\quad q_1= b_1x^{k_1}, \quad \dots \quad q_n= b_nx^{k_n}, \qquad k_0,k_1,\dots, k_n {\rm \;\; positive \; integers}.$$
Let us assume that the algorithm does not stop and goes on with an  infinite sequence of cases {\bf 4.2.3}; then the Puiseux series
$$p(x)= \sum_{n\geq 0} b_nx^{k_0+\dots+k_n}$$
is in fact a power series giving a unique 1-branch counted 3 times: but our curve $\bC$ is reduced and this is not possible. 
\dd Hence the only possible way in which {\bf 4.2.3} can appear is as follows: we have a finite sequence of $n$ steps {\bf 4.2.3} and then one of the other cases. 
\a If the step  $n+1$ is a case {\bf 4.1},  $\bC$ has one 3-branch at $O$; we have found:
 $$q_0= b_0x^{k_0},\; q_1= b_1x^{k_1}, \; \dots \; q_n= b_nx^{k_n}, \; q_{n+1}= c_1x^{i_1\over 3},\quad k_0,k_1,\dots, k_n \geq 0 {\rm  \; integers}, i_1\not\equiv 0\Mod{3}$$
which gives us $$p(x)=\sum_{\ell=0}^n b_{\ell}x^{k_0+\dots k_{\ell}}+c_1x^{k_0+\dots+k_n+{i_1\over 3}}$$
and since the common lowest denominator of $p(x)$ is 3, a parameterization for the 3-branch is 
$$\left(T^3, b_0T^{3k_0}+\dots+b_nT^{3(k_0+ \dots+k_n)}+ c_1T^{3(k_0+ \dots+k_n)+{i_1}}+ \dots \right)$$
with $ k_0,k_1,\dots, k_n \geq 0$ integers, and $ i_1\not\equiv 0\Mod{3}$. So we get $(i)$. 

\a  If the step  $n+1$ is not case {\bf 4.1}, then it has to be one of the cases from {\bf 4.2.1} to {\bf 7} (except {\bf4.2.3}), since the height of the Newton Polygon is still 3. It is easily seen that cases {\bf 4.2.1, 5.2.1, 6.2.1, 5.1, 6.1} give $(ii)$ and $(iii)$ just following the discussion in Remark  \ref{NP cases}.
\a  If the step $n+1$ is one of the cases {\bf 4.2.2, 5.2.2, 6.2.2}, in order to decide if $\bC$ has one 1-branch and one 2-branch or three 1-branches at $O$, the sub-algorithm relative to the double root must go on; now the height of the Newton Polygon is decreased to 2, so that we have to use cases from {\bf 2.1} to {\bf 3}, and again it is immediate to see that $(iv)$ holds following the discussion in the previous paragraph \ref{NP cases}.

\a Notice that we never use case {\rm \bf 1}, unless we wish to go on with the approximation of one of the 1-branches.
\prfend

\begin{defn}\label{type} In case $(i)$ we say that the triple point $O$ is of type $s$. This definition makes sense: if  $(T^3,p(T)), (T^3,\tilde p(T))$ are two equivalent parameterization of the 3-branch, by remark \ref{equiv}  we have the same value of $s$ for the two parameterizations.
\end{defn}

\medskip
We point out this definition because it is not hard to show that in case $(i)$ the branch is analytically equivalent to one of type $\left(T^3, c_1T^{s}+ \dots \right)$; it is sufficient to consider the biolomorphy of the plane around $O$ given by $\phi (x,y) =(x, -(b_0x^{h_0}+\dots+b_nx^{h_n}) + y)$, which 
maps the $3$-branch $(i)$ onto a $3$-branch of the desired form. Actually, more can be said about analytical equivalence in this case, see \cite{HH}, page 4.

\a {\small Stefano Canino, Dipartimento di Scienze Matematiche - Politecnico di Torino, stefano.canino@polito.it

\noindent Alessandro Gimigliano, Dipartimento di Matematica - Università di Bologna, alessandr.gimigliano@unibo.it

\noindent Monica Idà, Dipartimento di Matematica - Università di Bologna, monica.ida@unibo.it}

\end{document}